\newif\ifreview
\reviewfalse

\ifreview
  \documentclass[review]{siamart251104}
\else
  \documentclass[final]{siamart251104}
\fi
\newif\ifacm
\newif\ifarxiv
\arxivtrue

\usepackage[utf8]{inputenc}
\usepackage[english]{babel}

\ifacm
\usepackage{balance} 
\else
\usepackage{amsmath}
\usepackage{amsfonts}
\usepackage{amssymb}
\usepackage{amsopn}
\fi

\ifacm
\else
\fi

\ifacm
\else
\usepackage{booktabs}
\fi
\usepackage{tabularx}

\usepackage{color}
\usepackage{bm}
\usepackage{lineno}
\usepackage{url}
\usepackage{hyperref}
\usepackage{cleveref}
\usepackage[boxruled,algo2e]{algorithm2e}
\usepackage{setspace}
\usepackage{xcolor}
\usepackage{wrapfig}
\usepackage{graphicx}
\usepackage{pstool}
\usepackage{multirow}
\usepackage{setspace}
\usepackage{enumitem}
\usepackage{multirow}
\usepackage{subcaption}
\usepackage{comment}
\usepackage{xspace} 
\usepackage{tablefootnote} 
\usepackage{soul}   

\graphicspath{{./figures/}}

\DeclareCaptionFormat{custom}
{%
    \textbf{#1#2}\textnormal{\small #3}
}
\usepackage{pgfplots}
\usepgfplotslibrary{external} 
\pgfplotsset{compat=newest}

\usepackage{tikz}
\usetikzlibrary{shapes,snakes}
\usetikzlibrary{patterns}
\usetikzlibrary{shapes.misc}
\usetikzlibrary{arrows}
\usepgfplotslibrary{groupplots}
\usepgfplotslibrary{dateplot}
\usetikzlibrary{backgrounds}
\usetikzlibrary{decorations.pathreplacing,calligraphy}
\usetikzlibrary{external}

\newcommand{\code}[1]{\texttt{#1}}

\newcommand{\e}{\mathbf{e}}

\newcommand{\bk}{\mathbf{k}}

\newcommand{\flups}[1]{\code{FLUPS}#1}

\newcommand{\x}{\mathbf{x}}
\newcommand{\bx}{\mathbf{x}}

\newcolumntype{Y}{>{\centering\arraybackslash}X}
\newcolumntype{Z}[1]{>{\hsize=#1\arraybackslash}X}

\newcommand{\be}{\begin{equation}}
\newcommand{\eec}{\;\;,\end{equation}}
\newcommand{\ee}{ \end{equation}}
\newcommand{\eed}{\;\;.\end{equation}}
\newcommand{\bse}{\begin{subequations}}
\newcommand{\ese}{\end{subequations}}

\newcommand{\Tref}[1]{\Cref{#1}}
\newcommand{\Fref}[1]{\Cref{#1}}
\newcommand{\Eqqref}[1]{\Cref{#1}}

\newcommand{\Sref}[1]{\Cref{#1}}

\newcommand{\tref}[1]{\cref{#1}}
\newcommand{\fref}[1]{\cref{#1}}
\newcommand{\eqqref}[1]{\cref{#1}}

\newcommand{\sref}[1]{\cref{#1}}

\definecolor{darkGreen}{RGB}{0,160,0}
\definecolor{darkBlue}{RGB}{51,51,204}
\definecolor{darkRed}{RGB}{160,0,0}
\definecolor{orange}{RGB}{255,120,0}
\definecolor{myGreen}{rgb}{0.0706, 0.5333, 0.1921}

\ifreview

    \newcommand{\delete}[1]{\textcolor{darkRed}{\st{#1}}}
\else

    \newcommand{\delete}[1]{}
\fi

\definecolor{airforceblue}{rgb}{0.36, 0.54, 0.66}
\definecolor{carolinablue}{rgb}{0.6, 0.73, 0.89}
\definecolor{moonstoneblue}{rgb}{0.45, 0.66, 0.76}
\definecolor{amber}{rgb}{1.0, 0.75, 0.0}
\definecolor{cadmiumorange}{rgb}{0.93, 0.53, 0.18}
\definecolor{persianorange}{rgb}{0.85, 0.56, 0.35}
\definecolor{cadetgrey}{rgb}{0.57, 0.64, 0.69}
\definecolor{darkgray}{rgb}{0.66, 0.66, 0.66}
\definecolor{asparagus}{rgb}{0.53, 0.66, 0.42}
\definecolor{cambridgeblue}{rgb}{0.64, 0.76, 0.68}
\definecolor{olivine}{rgb}{0.6, 0.73, 0.45}
\definecolor{carnelian}{rgb}{0.7, 0.11, 0.11}
\definecolor{cardinal}{rgb}{0.77, 0.12, 0.23}
\definecolor{firebrick}{rgb}{0.7, 0.13, 0.13}

\definecolor{matlabBlue}{rgb}{0, 0.4470, 0.7410}
\definecolor{matlabOrange}{rgb}{0.8500, 0.3250, 0.0980}
\definecolor{matlabYellow}{rgb}{0.9290, 0.6940, 0.1250}
\definecolor{matlabPurple}{rgb}{0.4940, 0.1840, 0.5560}
\definecolor{matlabGreen}{rgb}{0.4660, 0.6740, 0.1880}
\definecolor{matlabRed}{rgb} {0.7765, 0.1569, 0.0902}
\definecolor{matlabBlack}{rgb}{0, 0, 0}

\definecolor{keyBlue}{rgb}{0.0392, 0.2274, 0.4235}
\definecolor{keyBlueTrans}{rgb}{0.4274, 0.4980, 0.5764}
\definecolor{keyGreen}{rgb}{0.05490, 0.38039, 0.003921}
\definecolor{keyGreenTrans}{rgb}{0.7019, 0.8313, 0.7176}
\definecolor{keyOrange}{rgb}{0.98039, 0.4980, 0.03529}
\definecolor{keyGray}{rgb}{0.50196, 0.50196, 0.50196}
\definecolor{keyGreen2}{rgb}{0.0627 , 0.4196 , 0.3882}
\definecolor{keyOrange2}{rgb}{0.9843 , 0.2941 , 0.2431}

\definecolor{pyBlue}{rgb}{0.11372,0.42352,0.67059}
\definecolor{pyBlueLight}{rgb}{0.6196,0.79215,0.88235}
\definecolor{pyOrange}{rgb}{1.0,0.4549,0.0627}
\definecolor{pyOrangeLight}{rgb}{0.9921568,0.643137,0.376470}
\definecolor{pyGreen}{rgb}{0.1529,0.5882,0.1529}
\definecolor{pyRed}{rgb}{0.81569,0.13725,0.14117}
\definecolor{pyPurple}{rgb}{0.53725,0.36078,0.7098}
\definecolor{myCacaDoie}{rgb}{0.580392156862745,0.43,0.27}
\definecolor{pyBrown}{rgb}{0.549019607843137,0.337254901960784,0.294117647058824}
\definecolor{pyRose}{rgb}{0.890196078431372,0.466666666666667,0.76078431372549}
\definecolor{pyOrange2}{rgb}{0.9921,0.5098,0.2078}
\definecolor{pyBlue2}{rgb}{0.3765,0.6431,0.8157}
\definecolor{pyBar1}{rgb}{0.1647,0.2745,0.2745}
\definecolor{pyBar2}{rgb}{0.0000,0.4588,0.4588}
\definecolor{pyBar3}{rgb}{0.4863,0.7804,0.9998}
\definecolor{pyBar4}{rgb}{0.4549,0.4549,0.4549}

\tikzset{cross/.style={cross out, draw, minimum size=2*(#1-\pgflinewidth), inner sep=0pt, outer sep=0pt}}

\newcommand{\Lop}{\mathcal{L}}
\newcommand{\Rop}{\mathcal{R}}

\newcommand{\sigL}{\sigma_{\Lop}}
\newcommand{\sigR}{\sigma_{\Rop}}
\newcommand{\Ghat}{\widehat{G}}
\newcommand{\bn}{\mathbf{n}}
\newcommand{\bmm}{\mathbf{m}}
\newcommand{\bp}{\mathbf{p}}
\newcommand{\bj}{\mathbf{j}}
\newcommand{\eul}{\mathrm{e}}
\newcommand{\ii}{\mathrm{i}}
\newcommand{\dd}{\,\mathrm{d}}
\newcommand{\Zbb}{\mathbb{Z}}

\newcommand{\Linf}[1]{\left\| #1 \right\|_{\infty}}

\newcommand{\lgf}[1]{\code{LGF#1}}
\newcommand{\mehr}[1]{\code{MEHR#1}}

\newcommand{\hej}[1]{\code{HEJ#1}}

\definecolor{col_lgf2}{RGB}{31,119,180}
\definecolor{col_lgf4}{RGB}{255,127,14}
\definecolor{col_lgf6}{RGB}{44,160,44}
\definecolor{col_lgf8}{RGB}{214,39,40}
\definecolor{col_mehr4}{RGB}{148,103,189}
\definecolor{col_mehr6}{RGB}{140,86,75}
\definecolor{col_hej10}{RGB}{23,190,207}
\definecolor{col_chat2}{RGB}{127,127,127}

\pgfplotsset{
    lgfplot/.style={
        width=0.48\textwidth, height=0.42\textwidth,
        grid=major, grid style={gray!25},
        tick label style={font=\scriptsize},
        label style={font=\small},
        legend style={font=\tiny, anchor=north east, at={(0.98,0.98)}, cells={anchor=west}},
        every axis plot/.append style={thick, mark size=1.6pt},
    },
    orderline/.style={black, dashed, thin, forget plot},
}
\newcommand{\balpha}{\bm{\alpha}}
\newcommand{\bbeta}{\bm{\beta}}

\title{Lattice Green's functions for Poisson problems with two unbounded and one periodic direction}

\author{
    Gilles Poncelet \thanks{Institute of Mechanics, Materials and Civil Engineering, Université catholique de Louvain, Belgium (corresponding author: \email{gilles.poncelet@uclouvain.be}).}
    \and Philippe Chatelain \footnotemark[1]
}

\begin{document}

\maketitle

\begin{abstract}
    Fourier-based Poisson solvers treat periodic directions spectrally and unbounded ones by the Hockney--Eastwood convolution with a free-space kernel. With two unbounded and one periodic direction, the setting of vortex tubes, wakes and jets, the kernel used so far by these solvers combines the two-dimensional kernel on the zero periodic mode with the symbol of the operator sampled on the doubled grid on every other mode. We show that the second ingredient is not exact: the error is the field of the periodic images of the source and independent of the mesh size. While the period $L_z$ is not much longer than the unbounded extent $L_x$, it is of order $\exp(-2\pi d/L_z)$, with $d$ the distance to the nearest image; beyond, it grows like $(L_z/L_x)^2$. The error is thus small on cubic domains, but larger than the solution itself on elongated ones. We remove the approximation by extending the integral representation of the exact second-order kernel of Hou and Colonius: the screened kernel of each mode is evaluated in closed form along one unbounded direction and by Gauss--Legendre quadrature along the other, for the high-order dimension-split and Mehrstellen stencils, at a one-off cost of about one solve. The kernels, implemented in the open-source library \flups{}, are verified by comparison with independent quadratures, by a padding test sensitive to the boundary condition enforced, and by convergence at the order of each stencil on an axially modulated tube source, on domains up to $64$ times longer than wide and up to $2.7\times10^{11}$ unknowns.
\end{abstract}



\section{Introduction}
\label{section:intro}
We consider the Poisson equation
\begin{equation}
    -\nabla^2 u = f \label{eq:poisson}
\end{equation}
on a Cartesian lattice, for a source $f$ of compact support, when some directions are unbounded (the solution decays at infinity) and the others are periodic. This combination is the natural setting of vortex tubes, wakes and jets in incompressible flow simulations~\cite{rennich_numerical_1997, chatelain_fourier-based_2010, rasmussen_multiresolution_2011, juul_spietz_regularization_2018, balty_2025}. It is met again in electrostatics, where the Coulomb potential of point charges periodic in one direction is evaluated by FFT-based Ewald methods~\cite{saffar_shamshirgar_spectral_2017, saffar_shamshirgar_fast_2021} that transform the periodic direction and treat the unbounded ones with a truncated continuum kernel~\cite{vico2016191}. Electronic-structure codes meet the same wire and surface geometries (one and two periodic directions, respectively) and solve them, screened or not, with interpolating scaling functions~\cite{genovese_efficient_2007, cerioni_efficient_2012}.

Fourier-based solvers treat the periodic directions spectrally and the unbounded ones with the Hockney--Eastwood algorithm~\cite{hockney_computer_2021, eastwood_remarks_1979}: the domain is doubled, the source zero-padded, and the free-space Green's function, truncated and mirrored, is convolved with the source by fast Fourier transforms. Provided the Green's function is known on the doubled domain, the circular convolution returns the free-space one without error, and the accuracy of the solve is that of the kernel alone: with the exact Green's function of the discrete operator, the difference equation is solved to round-off~\cite{caprace_flups_2021}.

The mixed case was first addressed by Rennich and Lele~\cite{rennich_numerical_1997}, who corrected a fully periodic solve with an analytical far-field expansion in modified Bessel functions. It was brought into the Hockney--Eastwood framework by Chatelain and Koumoutsakos~\cite{chatelain_fourier-based_2010}, who transformed the periodic direction, recognised each of its modes as a two-dimensional modified Helmholtz problem, and convolved it with the exact continuous kernel $K_0(|k_z| r)/2\pi$ evaluated in real space. Two families of kernels have followed. The regularised kernels of Hejlesen et al.~\cite{hejlesen_high_2013} are continuous Green's functions smoothed at the grid scale; their extension to mixed unbounded--periodic domains by Spietz et al.~\cite{juul_spietz_regularization_2018} sets out to follow the same route, but, having found no closed form for the regularised counterpart of that kernel, substitutes the two-dimensional kernel on the mean mode and the regularised symbol sampled at the frequencies of the doubled grid on every other mode.

The lattice Green's functions (LGFs) are the exact inverses of the finite-difference Laplacian on the infinite lattice~\cite{gillman_fast_2010, liska2014, liska_fast_2016}; their convolution with the source solves the difference equation to machine precision, which is what a consistent coupling with a finite-difference discretisation requires~\cite{gabbard_lattice_2023, dorschner_fast_2020, yu_multi-resolution_2022, gabbard_high-order_2024}. Gabbard and van Rees~\cite{gabbard_lattice_2023} recently derived the LGFs of the high-order dimension-split stencils and of the compact Mehrstellen operators~\cite{collatz1960numerical, deriaz_compact_2020} on fully unbounded domains, by fast asymptotic expansions, and on domains with one unbounded direction, in closed form by contour integration. They explicitly leave open the case of two unbounded and one periodic direction, ``which appear in the study of vortex tubes and the wakes of cylinders'', and note that the sampled-symbol approximation of~\cite{juul_spietz_regularization_2018} is one they did not pursue. That case was subsequently settled for the second-order stencil by Hou and Colonius~\cite{hou_fast_2025, hou_adaptive_2024}: each non-zero mode of the periodic direction is a two-dimensional discrete screened problem, whose exact LGF they obtain from an asymptotic expansion when the screening is strong and, otherwise, from a one-dimensional integral representation to which the trapezoidal rule converges exponentially, evaluated in batches by a fast Fourier transform. That kernel drives an adaptive LGF solver for flows past cylinders~\cite{hou_adaptive_2024}. Their construction covers the seven-point stencil, the second-order member of the family considered here, and their solver applies the kernel by a lattice fast multipole method rather than by FFT. What is missing is therefore an exact kernel for the high-order and compact operators on two unbounded and one periodic direction, usable in a Hockney--Eastwood solver, for which only the sampled-symbol approximation of~\cite{juul_spietz_regularization_2018} exists.

This paper closes the remaining gap in two steps. By the Poisson summation formula, sampling the symbol of a screened operator at the frequencies of the doubled grid amounts to convolving with the spatially periodised Green's function of that operator. We show that the approximation therefore commits an error equal to the field of the periodic images of the source. Each non-zero periodic mode is governed by a screened Poisson (modified Helmholtz) operator whose Green's function decays exponentially, so the error decays exponentially with the distance to the nearest image. It is independent of the mesh size, small on a cubic domain (though above the discretisation error of the high-order kernels), and larger than the solution itself when the periodic direction is much longer than the unbounded ones, where the images add up.

We then remove the approximation altogether, by extending the construction of~\cite{hou_fast_2025}. On each mode, the screened LGF is obtained in closed form along one unbounded direction, from the one-unbounded LGFs of~\cite{gabbard_lattice_2023}, and to machine precision along the second unbounded direction by Gauss--Legendre quadrature on a half period, before the standard Hockney--Eastwood assembly. The construction requires no new special function and covers the dimension-split stencils of order two to eight, the Mehrstellen operators of order four and six and, with the regularised one-unbounded kernel of~\cite{juul_spietz_regularization_2018} as integrand, the regularised kernels of order two to ten; of these only the second-order LGF was available before~\cite{hou_fast_2025}. The kernel is set up once per grid, at a cost of about one solve at most on the grids considered here.

The paper is organised as follows. \Sref{section:method} recalls the lattice setting, derives the error of the sampled-symbol approximation and describes the exact construction. \Sref{section:validation} verifies the kernels against independent quadratures, against a padding test that is sensitive to the boundary condition actually enforced, and on an axially modulated tube source whose solution is known semi-analytically, on domains up to $64$ times longer than wide and up to $4096^2\times16384$ unknowns. \Sref{section:conclusion} concludes; Appendix~\ref{app:numerics} collects the numerical details: the evaluation of the closed forms, the quadrature rules and the reference solution of the tube test.


\section{Methodology}
\label{section:method}

\subsection{Lattice Green's functions and the Hockney--Eastwood algorithm}
\label{sec:lgf}
We follow the notation of~\cite{gabbard_lattice_2023}. On a lattice of spacing $h$ with $N_x\times N_y\times N_z$ points, a finite-difference discretisation of \eqqref{eq:poisson} reads
\begin{equation}
    \left[\Lop u\right](\bn) = \left[\Rop f\right](\bn), \qquad
    \left[\Lop u\right](\bn) = \frac{1}{h^2}\sum_{\balpha} a_{\balpha}\, u(\bn+\balpha),\quad
    \left[\Rop f\right](\bn) = \sum_{\bbeta} b_{\bbeta}\, f(\bn+\bbeta),
    \label{eq:stencil}
\end{equation}
with integer multi-indices $\balpha,\bbeta$. For the classical dimension-split stencils of order $2,4,6,8$ (\lgf{2}--\lgf{8}) $\Rop$ is the identity; for the compact Mehrstellen operators of order four and six~\cite{collatz1960numerical, spotz_high-order_1996, sutmann_high-order_2006, deriaz_compact_2020} (\mehr{4}, \mehr{6}) $\Lop$ is a 19- or 27-point stencil of radius one and $\Rop$ a small averaging operator. Both operators are symmetric and consistent, so their Fourier symbols
\begin{equation}
    \sigL(\bk) = \sum_{\balpha} a_{\balpha}\,\eul^{\ii\bk\cdot\balpha}, \qquad
    \sigR(\bk) = \sum_{\bbeta} b_{\bbeta}\,\eul^{\ii\bk\cdot\bbeta}, \qquad \bk\in[-\pi,\pi]^3,
\end{equation}
are real and even, $\sigL(\bk) = |\bk|^2 + \mathcal{O}(|\bk|^4)$ and $\sigL(\bk)>0$ for $\bk\neq 0$. Both families are built from the symbol of minus the centred second difference, $s^2(k) = 4\sin^2(k/2)$; we introduce the notation $s_i^2 = s^2(k_i)$ with $i\in\{x,y,z\}$ to make its dimension explicit. The dimension-split symbols are then sums $\sigL = P(s_x^2)+P(s_y^2)+P(s_z^2)$ of the symbol of the one-dimensional stencil, $P(s^2) = s^2 + s^4/12 + s^6/90 + s^8/560$ truncated after one to four terms for \lgf{2}--\lgf{8}. The Mehrstellen symbols are~\cite{gabbard_lattice_2023}
\begin{equation}
    \begin{aligned}
        \sigL^{(4)} &= (s_x^2+s_y^2+s_z^2) - \tfrac{1}{6}\,(s_x^2s_y^2+s_x^2s_z^2+s_y^2s_z^2), \qquad
        \sigL^{(6)} = \sigL^{(4)} + \tfrac{1}{30}\, s_x^2s_y^2s_z^2, \\
        \sigR^{(4)} &= 1 - \tfrac{1}{12}\,(s_x^2+s_y^2+s_z^2), \\
        \sigR^{(6)} &= \sigR^{(4)}  + \tfrac{1}{90}\,(s_x^2s_y^2+s_x^2s_z^2+s_y^2s_z^2) - \tfrac{1}{240}\,(s_x^4+s_y^4+s_z^4).
    \end{aligned}
    \label{eq:mehr_symbols}
\end{equation}
The resolution of the Mehrstellen-discretized system of \eqqref{eq:stencil} can be performed in two fashions. One may obtain the so-called Laplacian-only kernel (\mehr{4L}, \mehr{6L}) through the inversion of $\Lop$ alone; the expected input for the solver is then $f^\star = \Rop f$. Conversely, the folded kernel (\mehr{4F}, \mehr{6F}) carries $\sigR$ in its numerator and returns $\Lop^{-1}\Rop f$ from the raw $f$.

The lattice Green's function (LGF) is the decaying solution of $\Lop G = \Rop\delta$,
\begin{equation}
    G(\bn) = \frac{h^2}{(2\pi)^3}\int_{[-\pi,\pi]^3} \eul^{\ii\bn\cdot\bk}\,\frac{\sigR(\bk)}{\sigL(\bk)}\dd\bk,
    \label{eq:lgf3d}
\end{equation}
and the solution of \eqqref{eq:stencil} on the infinite lattice is the convolution $u = G * f$. On a domain that is unbounded in some directions and periodic in the others, the convolution is evaluated by the Hockney--Eastwood (HE) algorithm~\cite{hockney_computer_2021}: the domain is doubled in each unbounded direction, the source is zero-padded, and the kernel is truncated to $0\le n\le N$, mirrored about $n = N$ and transformed. The circular convolution then equals the free-space one on the original box: the kernel being truncated, no periodic image of the source contributes there, i.e.\ the convolution is free of aliasing. We refer to this justification as the HE argument. The periodic directions are transformed directly and no padding is needed there. All that is required is the LGF on the doubled domain, fully transformed in the periodic directions.

The fully unbounded and the one-unbounded cases are covered by~\cite{gabbard_lattice_2023}, by the asymptotic expansions of Martinsson and Rodin~\cite{martinsson_asymptotic_2002} and by closed forms that generalise Buneman's~\cite{buneman_analytic_1971}. The intermediate case of two unbounded and one periodic direction is the subject of this paper; for the second-order stencil it has been settled by Hou and Colonius~\cite{hou_fast_2025, hou_adaptive_2024}, on lattices that may be anisotropic. Their integral representation~\cite[Thm.~4.1]{hou_fast_2025} is the seven-point instance of the construction of \sref{sec:exact}, taken over the full period rather than a half period. Its derivation, a discrete Fourier transform along one unbounded direction and a recurrence along the other, applies to any stencil whose one-unbounded LGF is known in closed form; \sref{sec:exact} uses those of~\cite{gabbard_lattice_2023}. The classical closed forms of the same object in terms of hypergeometric functions go back to Katsura and Inawashiro~\cite{katsura_lattice_1971}.

\subsection{Slice decomposition}
\label{sec:slices}
Let $z$ be the periodic direction, of period $L_z = N_z h$. Since the symbol is a polynomial in $s^2_z$, the discrete Fourier transform along $z$ alone block-diagonalises the problem: the mode $k_z = 2\pi j_z/N_z$ obeys a two-dimensional lattice equation in the $(x,y)$ plane whose symbol is $\sigL(k_x,k_y,k_z)$ at fixed $k_z$, and whose Green's function is the \emph{slice} LGF
\begin{equation}
    G_{k_z}(n_x,n_y) = \frac{h^2}{(2\pi)^2}\iint_{[-\pi,\pi]^2} \eul^{\ii(n_xk_x+n_yk_y)}\,
    \frac{\sigR(k_x,k_y,k_z)}{\sigL(k_x,k_y,k_z)}\dd k_x\dd k_y .
    \label{eq:slice}
\end{equation}
For $k_z\neq0$ the slice symbol has a positive minimum, since $\sigL>0$ for $\bk\neq0$. It no longer vanishes at $k_x = k_y = 0$, where it takes the value $\sigL(0,0,k_z)\simeq k_z^2$: each slice is a discrete \emph{screened} Poisson (modified Helmholtz) problem. The integrand of \eqqref{eq:slice} is then analytic and $G_{k_z}$ decays exponentially, like the continuous screened kernel $\frac{1}{2\pi}K_0(\kappa|\bn|)$:
\begin{equation}
    |G_{k_z}(\bn)| \lesssim \frac{\eul^{-\kappa|\bn|}}{\sqrt{\kappa|\bn|}}\quad (\kappa|\bn|\gtrsim1), \qquad
    \kappa = |k_z| + \mathcal{O}(k_z^3).
    \label{eq:decay}
\end{equation}
The decay rate $\kappa$ has this expansion for every consistent stencil, because $\sigL = |\bk|^2+\mathcal{O}(|\bk|^4)$: on the low modes, the screening length $\kappa^{-1}\simeq N_z/(2\pi j_z)$ does not depend on the kernel.

The mode $k_z = 0$ is the exception: the slice symbol vanishes quadratically at $k_x = k_y = 0$, and the slice LGF is the unscreened two-dimensional LGF of $\sigR/\sigL$, which grows logarithmically instead of decaying. This dichotomy between the mode $k_z = 0$ and the $N_z-1$ screened ones is why the kernel of Spietz et al.~\cite{juul_spietz_regularization_2018} has two cases. Singly periodic Ewald methods~\cite{saffar_shamshirgar_spectral_2017, saffar_shamshirgar_fast_2021} set the mode $k_z = 0$ apart for the same reason: they treat it on its own, and need less upsampling of the unbounded directions for the screened modes.

\subsection{The sampled-symbol approximation and its error}
\label{sec:error}
The approximation of~\cite{juul_spietz_regularization_2018}, transposed to a lattice symbol, defines the spectral kernel on the doubled grid, hereafter the \emph{Spietz kernel}, as
\begin{equation}
    \Ghat(k_x,k_y,k_z) =
    \begin{cases}
        \Ghat_{\mathrm{2D}}(k_x,k_y), & k_z = 0,\\[2pt]
        h^2\,\sigR(k_x,k_y,k_z)/\sigL(k_x,k_y,k_z), & k_z\neq0,
    \end{cases}
    \qquad k_x = \frac{\pi j_x}{N_x},\ k_y = \frac{\pi j_y}{N_y},
    \label{eq:formula}
\end{equation}
where $\Ghat_{\mathrm{2D}}$ is the two-dimensional LGF, tabulated near the origin and expanded in the far field, truncated, mirrored and transformed on the doubled grid, and where the $k_z \neq 0$ case is the symbol \emph{sampled} at the frequencies of the $2N_x\times2N_y$ grid. The first case is exact. The second one is not, and it is presented in~\cite{juul_spietz_regularization_2018} as a deliberate fallback, the closed form of the mixed kernel having proved elusive: ``it is however not trivial to find a closed form expression of the integral [\dots], but we notice that to compute discrete solutions [\dots] using the FFT-based Poisson solver a regular spectral representation is sufficient. The integrand [\dots] is regular for [$k_z\neq0$] and may be used as it is stated.''

\begin{figure}[htbp]
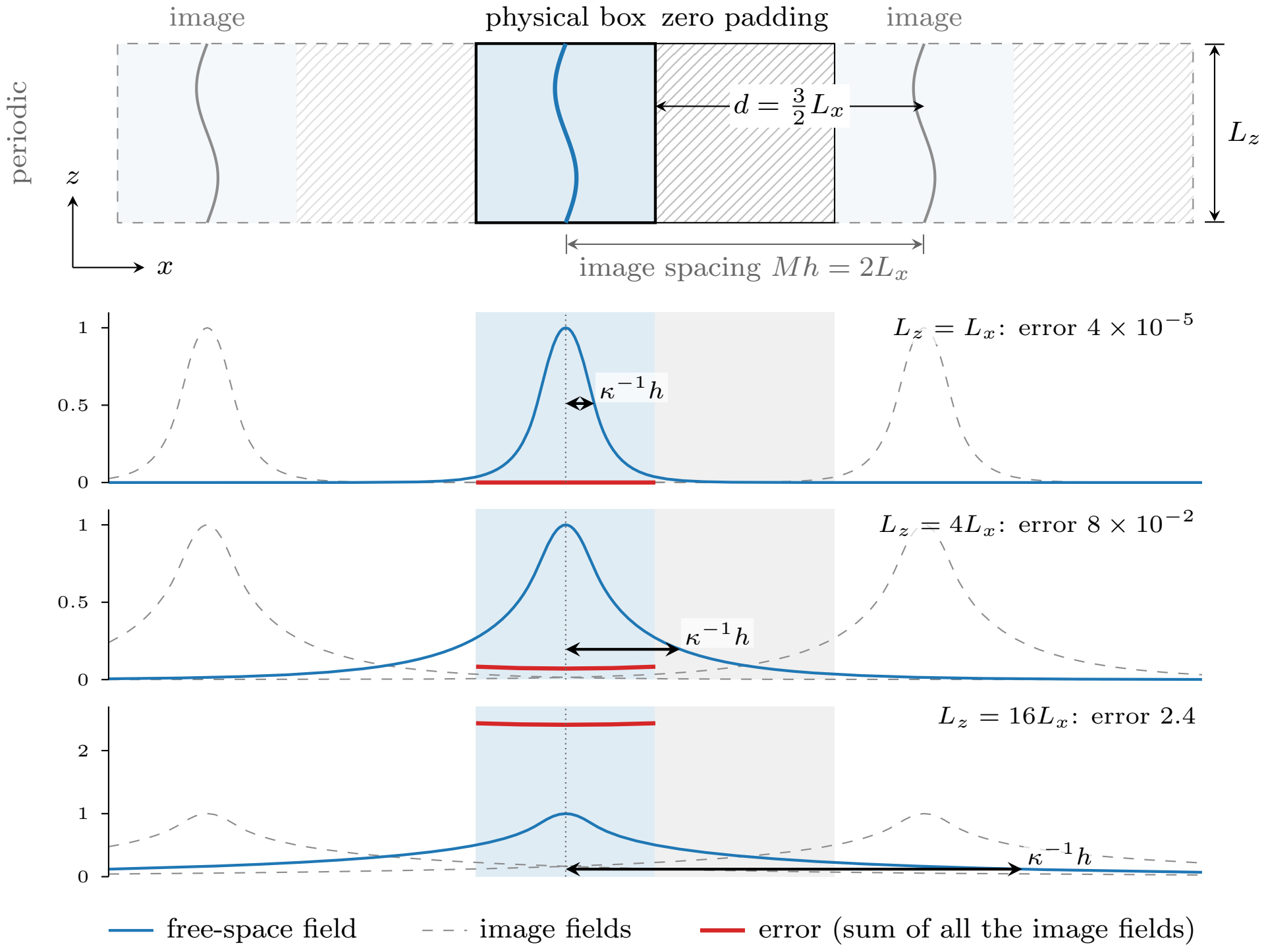

    \centering
\begin{tikzpicture}[x=1.8cm, y=1.8cm, >=stealth, font=\footnotesize]
    \def\H{1.0}
    \def\xc{0.5}
    \def\sh{0.95}
    \begin{scope}[shift={(0,3.55)}]
    \foreach \s in {-2,2} {
        \fill[col_lgf2!5] (\s,0) rectangle (\s+1,\H);
        \fill[pattern=north east lines, pattern color=black!12] (\s+1,0) rectangle (\s+2,\H);
        \draw[black!45, dashed] (\s,0) rectangle (\s+2,\H);
        \draw[black!45, thick] plot[domain=0:\H, samples=60, smooth] ({\s+\xc+0.06*sin(360*\x/\H)}, \x);
    }
    \node[black!55, above] at (-1.5,\H) {image};
    \node[black!55, above] at (2.5,\H) {image};
    \fill[col_lgf2!14] (0,0) rectangle (1,\H);
    \fill[pattern=north east lines, pattern color=black!25] (1,0) rectangle (2,\H);
    \draw[thick] (0,0) rectangle (1,\H);
    \draw (1,0) rectangle (2,\H);
    \node[above] at (0.5,\H) {physical box};
    \node[above, align=center] at (1.5,\H) {zero padding};
    \draw[col_lgf2, very thick] plot[domain=0:\H, samples=60, smooth] ({\xc+0.06*sin(360*\x/\H)}, \x);
    \draw[<->] (1,0.65) -- (\xc+2,0.65) node[midway, fill=white, inner sep=1pt] {$d = \tfrac32L_x$};
    \draw[|<->|, black!60] (\xc,-0.12) -- (\xc+2,-0.12) node[midway, below] {image spacing $Mh = 2L_x$};
    \draw[|<->|] (4.12,0) -- (4.12,\H) node[midway, right] {$L_z$};
    \draw[->] (-2.25,-0.25) -- (-1.85,-0.25) node[right] {$x$};
    \draw[->] (-2.25,-0.25) -- (-2.25,0.15) node[above] {$z$};
    \node[rotate=90, black!60] at (-2.55,0.5) {periodic};
    \end{scope}
    \newcommand{\imgstrip}[8]{
        \begin{scope}[shift={(0,#1)}]
            \fill[col_lgf2!14] (0,0) rectangle (1,\sh);
            \fill[black!6] (1,0) rectangle (2,\sh);
            \begin{scope}
                \clip (-2.05,-0.03) rectangle (4.05,\sh);
                \foreach \p in {m1,p1}
                    {\input{data/images_sketch/Lz#2_img\p.tex}\draw[black!45, dashed] \sketchcurve;}
                \input{data/images_sketch/Lz#2_free.tex}\draw[col_lgf2, thick] \sketchcurve;
                \input{data/images_sketch/Lz#2_err.tex}\draw[col_lgf8, very thick] \sketchcurve;
            \end{scope}
            \draw[densely dotted, black!60] (\xc,0) -- (\xc,\sh);
            \draw[<->, thick] (\xc,#4) -- ({\xc+#3},#4) node[above right, font=\scriptsize, inner sep=1pt, fill=white, fill opacity=0.7, text opacity=1] {$\kappa^{-1}h$};
            \draw (-2.05,0) -- (-2.05,\sh);
            \foreach \e in {#7}
                \draw (-2.05,{\e/#6*\sh}) -- ++(-0.04,0) node[left, font=\tiny] {$\e$};
            \node[#8, font=\scriptsize, align=right, fill=white, fill opacity=0.7, text opacity=1, inner sep=1.5pt] {#5};
        \end{scope}
    }
    
        \begin{scope}[shift={(0,2.1)}]
            \fill[col_lgf2!14] (0,0) rectangle (1,\sh);
            \fill[black!6] (1,0) rectangle (2,\sh);
            \begin{scope}
                \clip (-2.05,-0.03) rectangle (4.05,\sh);
                \foreach \p in {m1,p1}
                    {\input{data/images_sketch/Lz1_img\p.tex}\draw[black!45, dashed] \sketchcurve;}
                \input{data/images_sketch/Lz1_free.tex}\draw[col_lgf2, thick] \sketchcurve;
                \input{data/images_sketch/Lz1_err.tex}\draw[col_lgf8, very thick] \sketchcurve;
            \end{scope}
            \draw[densely dotted, black!60] (\xc,0) -- (\xc,\sh);
            \draw[<->, thick] (\xc,0.441) -- ({\xc+0.159},0.441) node[above right, font=\scriptsize, inner sep=1pt, fill=white, fill opacity=0.7, text opacity=1] {$\kappa^{-1}h$};
            \draw (-2.05,0) -- (-2.05,\sh);
            \foreach \e in {0, 0.5, 1}
                \draw (-2.05,{\e/1.1*\sh}) -- ++(-0.04,0) node[left, font=\tiny] {$\e$};
            \node[anchor=north east, at={(4.05,\sh+0.03)}, font=\scriptsize, align=right, fill=white, fill opacity=0.7, text opacity=1, inner sep=1.5pt] {$L_z = L_x$: error $4\times10^{-5}$};
        \end{scope}

        \begin{scope}[shift={(0,1.0)}]
            \fill[col_lgf2!14] (0,0) rectangle (1,\sh);
            \fill[black!6] (1,0) rectangle (2,\sh);
            \begin{scope}
                \clip (-2.05,-0.03) rectangle (4.05,\sh);
                \foreach \p in {m1,p1}
                    {\input{data/images_sketch/Lz4_img\p.tex}\draw[black!45, dashed] \sketchcurve;}
                \input{data/images_sketch/Lz4_free.tex}\draw[col_lgf2, thick] \sketchcurve;
                \input{data/images_sketch/Lz4_err.tex}\draw[col_lgf8, very thick] \sketchcurve;
            \end{scope}
            \draw[densely dotted, black!60] (\xc,0) -- (\xc,\sh);
            \draw[<->, thick] (\xc,0.169) -- ({\xc+0.637},0.169) node[above right, font=\scriptsize, inner sep=1pt, fill=white, fill opacity=0.7, text opacity=1] {$\kappa^{-1}h$};
            \draw (-2.05,0) -- (-2.05,\sh);
            \foreach \e in {0, 0.5, 1}
                \draw (-2.05,{\e/1.1*\sh}) -- ++(-0.04,0) node[left, font=\tiny] {$\e$};
            \node[anchor=north east, at={(4.05,\sh+0.03)}, font=\scriptsize, align=right, fill=white, fill opacity=0.7, text opacity=1, inner sep=1.5pt] {$L_z = 4L_x$: error $8\times10^{-2}$};
        \end{scope}

        \begin{scope}[shift={(0,-0.1)}]
            \fill[col_lgf2!14] (0,0) rectangle (1,\sh);
            \fill[black!6] (1,0) rectangle (2,\sh);
            \begin{scope}
                \clip (-2.05,-0.03) rectangle (4.05,\sh);
                \foreach \p in {m1,p1}
                    {\input{data/images_sketch/Lz16_img\p.tex}\draw[black!45, dashed] \sketchcurve;}
                \input{data/images_sketch/Lz16_free.tex}\draw[col_lgf2, thick] \sketchcurve;
                \input{data/images_sketch/Lz16_err.tex}\draw[col_lgf8, very thick] \sketchcurve;
            \end{scope}
            \draw[densely dotted, black!60] (\xc,0) -- (\xc,\sh);
            \draw[<->, thick] (\xc,0.042) -- ({\xc+2.546},0.042) node[above right, font=\scriptsize, inner sep=1pt, fill=white, fill opacity=0.7, text opacity=1] {$\kappa^{-1}h$};
            \draw (-2.05,0) -- (-2.05,\sh);
            \foreach \e in {0, 1, 2}
                \draw (-2.05,{\e/2.7*\sh}) -- ++(-0.04,0) node[left, font=\tiny] {$\e$};
            \node[anchor=north east, at={(4.05,\sh+0.03)}, font=\scriptsize, align=right, fill=white, fill opacity=0.7, text opacity=1, inner sep=1.5pt] {$L_z = 16L_x$: error $2.4$};
        \end{scope}
    
    \begin{scope}[shift={(-2.05,-0.4)}]
        \draw[col_lgf2, thick] (0,0) -- (0.25,0) node[right, black] {free-space field};
        \draw[black!45, dashed] (1.75,0) -- (2.0,0) node[right, black] {image fields};
        \draw[col_lgf8, very thick] (3.3,0) -- (3.55,0) node[right, black] {error (sum of all the image fields)};
    \end{scope}
\end{tikzpicture}
    \caption{Error of the Spietz kernel in the $(x,z)$ plane. Top: physical box, zero padding and two image cells of the period $2L_x$ imposed by the Spietz kernel, \eqqref{eq:poisson_sum}, with a compact source of mode $j_z = 1$ and its images (to scale for $L_z = L_x$). Below, from top to bottom: field of that mode along $x$ for $L_z = L_x$, $4L_x$ and $16L_x$, relative to its peak, the error growing with the period (note the scale of the last strip); the arrows mark the screening length $\kappa^{-1}h\simeq L_z/2\pi$.}
    \label{fig:images}
\end{figure}

Regularity of the symbol is not the same thing as correctness of the boundary condition, and the mechanism by which the two dissociate is instructive. Let $M = 2N$ be the doubled length (in each unbounded direction). The Poisson summation formula states that sampling the transform of an absolutely summable sequence at the $M$-point DFT frequencies returns the $M$-point DFT of its $M$-periodised version, which is why free-space kernels are sampled in real space~\cite{hockney_computer_2021, vico2016191}:
\begin{equation}
    h^2\frac{\sigR}{\sigL}\!\left(\tfrac{2\pi\bj}{M},k_z\right)
    = \sum_{\bn\in[0,M)^2} G^{\mathrm{per}}_{k_z}(\bn)\,\eul^{-2\pi\ii\,\bn\cdot\bj/M},
    \qquad
    G^{\mathrm{per}}_{k_z}(\bn) = \sum_{\bp\in\Zbb^2} G_{k_z}(\bn+M\bp),
    \label{eq:poisson_sum}
\end{equation}
The second case of \eqqref{eq:formula} is therefore nothing but the ordinary periodic solve on the doubled box. On a mode $k_z\neq0$, let $u_{k_z}(\bn)$ and $f_{k_z}(\bn)$ be that mode of $u$ and $f$, i.e.\ the two-dimensional fields left by the transform along $z$ alone, and $u^{\mathrm{per}}_{k_z} = G^{\mathrm{per}}_{k_z} * f_{k_z}$ the result of this solve, restricted to the original box. Splitting the periodised kernel,
\begin{equation}
    u^{\mathrm{per}}_{k_z}(\bn) = \underbrace{\sum_{\bmm} G_{k_z}(\bn-\bmm)\,f_{k_z}(\bmm)}_{\text{free-space solution}}
    \;+\; \underbrace{\sum_{\bmm}\sum_{\bp\neq0} G_{k_z}(\bn-\bmm+M\bp)\,f_{k_z}(\bmm)}_{\text{field of the periodic images}}.
    \label{eq:split}
\end{equation}
The first sum is the exact solution, by the HE argument; the second one is the error (\fref{fig:images}). Since $\bn,\bmm\in[0,N)^2$, every image displacement satisfies $|\bn-\bmm+M\bp|\ge N+1$: the distance $d$ from a target point to the nearest image of the source is at least $(N+1)h\simeq L_x$. As long as the screening length of the lowest mode, $\kappa^{-1} \simeq L_z/2\pi h$, is short compared with the image spacing $M$ ($\kappa M\gg1$, i.e.\ $L_z\ll4\pi L_x$), the nearest images dominate (\fref{fig:images}, $L_z = 4L_x$), and \eqqref{eq:decay} yields, on the least screened mode $j_z = 1$, which dominates the error,
\begin{equation}
    \Linf{u^{\mathrm{per}}_{k_z}-u_{k_z}} \lesssim \underbrace{\|f_{k_z}\|_1}_{h^2\sum_{\bmm}|f_{k_z}(\bmm)|}\;\underbrace{\frac{\eul^{-\kappa d/h}}{\sqrt{\kappa d/h}}}_{\text{nearest image}}, \qquad \frac{\kappa d}{h}\simeq\frac{2\pi d}{L_z} \qquad (\kappa M\gg1).
    \label{eq:bound}
\end{equation}
Conversely, when $\kappa M\lesssim1$, i.e.\ $L_z\gtrsim4\pi L_x$, the screening length exceeds the image spacing: the kernel barely decays from one image to the next, and the images add up to the mean of the periodised kernel over the doubled cell (\fref{fig:images}, $L_z = 16L_x$):
\begin{equation}
    \begin{aligned}
        \Linf{u^{\mathrm{per}}_{k_z}-u_{k_z}} &\simeq \frac{1}{M^2}\,\overbrace{\frac{\sigR}{\sigL}(0,0,k_z)}^{\simeq\,(L_z/2\pi h)^2}\,\overbrace{\Bigl|h^2\sum_{\bmm} f_{k_z}(\bmm)\Bigr|}^{\text{net source}}\\
        &\simeq \Bigl(\frac{L_z}{4\pi L_x}\Bigr)^2\Bigl|h^2\sum_{\bmm} f_{k_z}(\bmm)\Bigr| \qquad (\kappa M\lesssim1).
    \end{aligned}
    \label{eq:bound_long}
\end{equation}
The error grows algebraically with $L_z$ and depends on the kernel only through $\sigL/\sigR(0,0,k_z) = k_z^2 + \mathcal{O}(k_z^{p+2})$, $p$ being the order of the kernel, i.e.\ not at all for small $k_z$. If the net source, i.e.\ the monopole, vanishes, as for a pair of opposite-signed compact sources, the next term, i.e.\ the dipole, takes over: the kernel of the image field in \eqqref{eq:split}, $\sum_{\bp\neq0}G_{k_z}(\bn+M\bp)$, is even and screened near the origin, so that it tends to a constant plus $h^2|\bx_{\bn}|^2/(16L_x^2)$, and the error levels off at
\begin{equation}
    u^{\mathrm{per}}_{k_z}(\bn)-u_{k_z}(\bn) \simeq \left(\frac{h}{4L_x}\right)^2\sum_{\bmm} f_{k_z}(\bmm)\,|\bx_{\bn}-\bx_{\bmm}|^2
\end{equation}
which no longer depends on $L_z$ or the kernel, and is about $|\boldsymbol{D}|/(16L_x)$ for a source centered in the box, $\boldsymbol{D} = h^2\sum_{\bmm} f_{k_z}(\bmm)\,\bx_{\bmm}$ being its dipole moment. Three consequences follow:
\begin{enumerate}[label=(\roman*), nosep, leftmargin=*]
    \item The approximation is not exact; at best, it is exponentially accurate when the lowest mode is screened.
    \item The error depends only on the geometry of the domain, not on $h$ or $N$: it therefore constitutes a lower bound for the error that is independent of any grid refinement.
    \item The controlling parameter is the aspect ratio. For a compact unit source centred in a cubic domain, the nearest image lies at $d\simeq\tfrac32L_x$ and the image field is $4\times10^{-5}$ of the solution; it is $8\times10^{-2}$ for $L_z = 4L_x$ and twice the solution for $L_z = 16L_x$ (\fref{fig:images}).
\end{enumerate}

The test invoked in support of \eqqref{eq:formula}, in~\cite{juul_spietz_regularization_2018} as in the validation of \flups{}~\cite{caprace_flups_2021, balty_flups_2023}, is a compact manufactured solution, and it is structurally blind to \eqqref{eq:split}. A manufactured solution that vanishes at the edges of the box is the solution of both the unbounded and the periodic problem, so the Spietz kernel and the exact one both return it to machine precision and the test shows clean convergence whatever the kernel. The discrete residual, which vanishes for any exact Green's function of the doubled periodic box, cannot separate the two either. A discriminating test must use a compact \emph{source} whose solution is $\mathcal{O}(1)$ on the boundary, and a periodic direction longer than the unbounded ones. \Sref{section:validation} is built on this observation.

\subsection{Exact screened slices}
\label{sec:exact}
The defect is confined to the screened slices, on which the approximation \eqqref{eq:formula} replaces the true kernel by its periodised version. It disappears if $G_{k_z}(n_x,n_y)$ is evaluated in physical space for $0\le n_x\le N_x$, $0\le n_y\le N_y$ and every $k_z\neq0$, and the classical HE recipe is then applied slice by slice: truncation to the box, mirroring about $N_x$ and $N_y$, and a two-dimensional transform, exactly as for the $k_z = 0$ slice. Nothing changes in the layout of the kernel, in the transforms or in the online solve. What is needed is a way to evaluate \eqqref{eq:slice} to machine precision for an arbitrary stencil; the following steps do so.

\paragraph{Step 1: the $k_x$ integral in closed form}
Freeze $(k_y,k_z)$. The inner integral
\begin{equation}
    I(n;k_y,k_z) = \frac{1}{2\pi}\int_{-\pi}^{\pi}\eul^{\ii n k}\,\frac{\sigR(k;k_y,k_z)}{\sigL(k;k_y,k_z)}\dd k
    \label{eq:inner}
\end{equation}
is the one-unbounded LGF of~\cite{gabbard_lattice_2023}, with $k_y$ and $k_z$ as its two spectral wavenumbers. Its integrand is a rational function of $\eul^{\ii k}$, with as many poles inside the unit circle as the degree of $\sigL$ in $\cos k$, so the integral has a closed form. That reference gives numerically stable ones for every stencil considered here, including the repeated-root and small-wavenumber regimes, and we use them. On elongated grids their evaluation calls for a few precautions against cancellation, given in Appendix~\ref{app:numerics}. With them, the one-dimensional LGF agrees with an extended-precision evaluation to $10^{-14}$ over the whole $(k_y,k_z,n)$ range, up to $n = 1024$ and $k_z = 2\pi/4096$.

The structure of these closed forms is what the next step relies on, and the stencils of radius one (\lgf{2}, \mehr{4}, \mehr{6}) display it in full. Their symbol is linear in $\cos k_x$,
\begin{equation}
    \sigL(k_x;k_y,k_z) = e_0(k_y,k_z) - e_1(k_y,k_z)\cos k_x,
    \label{eq:affine}
\end{equation}
with, for the seven-point stencil, $e_0 = 2 + s^2_y + s^2_z$ and $e_1 = 2$. Ellipticity of the discrete Laplace operator gives ($\sigL>0$ for $\bk\neq0$) $e_0 > |e_1|$ for $(k_y,k_z)\neq(0,0)$, and the closed form reduces, for $\Rop = I$ and $n\ge0$, to
\begin{equation}
    I(n;k_y,k_z) = \frac{\rho^{\,n}}{\sqrt{e_0^2-e_1^2}}, \qquad
    \rho = \frac{e_1}{e_0 + \sqrt{e_0^2-e_1^2}}, \qquad |\rho|<1 .
    \label{eq:residue}
\end{equation}
The inner integral thus decays geometrically in $n$, at a rate set by the distance of the slice symbol to zero, and it is analytic in $(k_y,k_z)$ wherever $e_0^2 > e_1^2$. The higher-order stencils share both properties, with one geometric term per root.

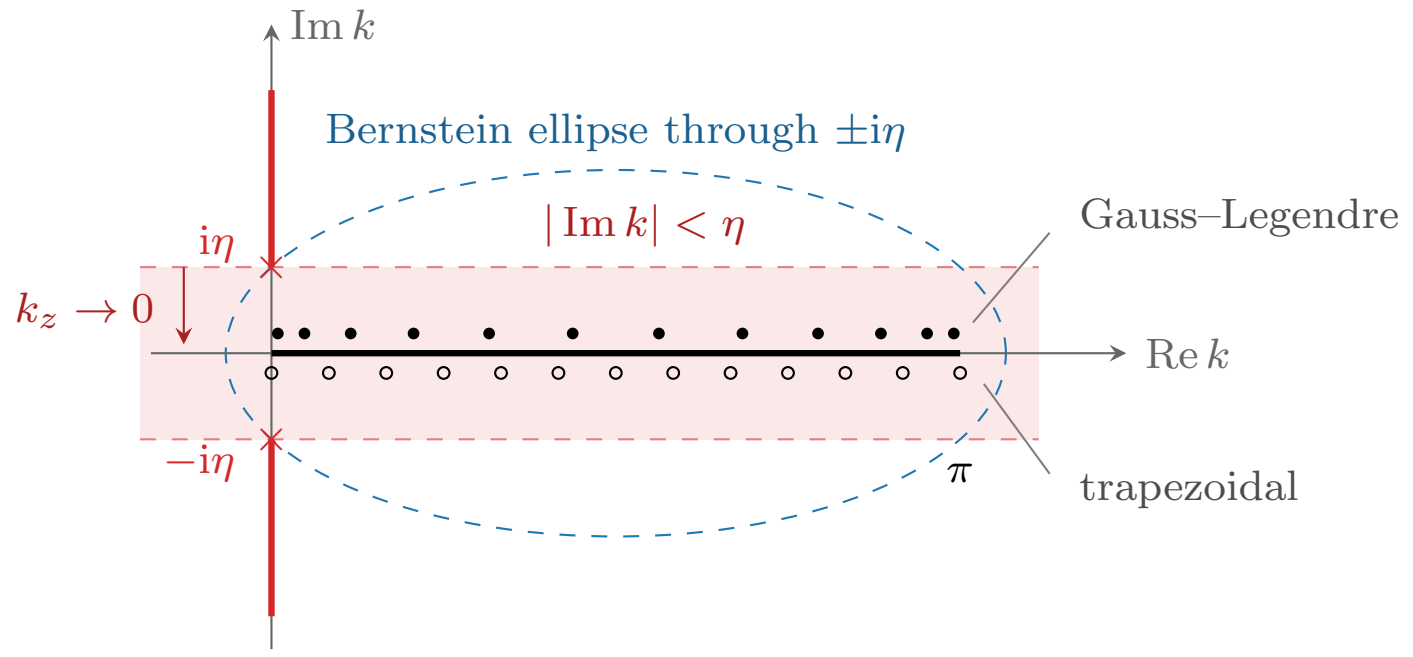
\begin{figure}[htbp]
    \centering
\begin{tikzpicture}[>=stealth, font=\footnotesize, x=1.45cm, y=1.45cm]
    \def\cut{1.2}
    \def\bp{0.3927}
    \fill[col_lgf8!10] (-0.6,-\bp) rectangle (3.5,\bp);
    \draw[col_lgf8!60, dashed] (-0.6,\bp) -- (3.5,\bp);
    \draw[col_lgf8!60, dashed] (-0.6,-\bp) -- (3.5,-\bp);
    \node[col_lgf8!80!black, above] at (1.7,0.40) {$|\operatorname{Im} k| < \eta$};
    \pgfmathsetmacro{\ea}{(\bp + sqrt(pi*pi + \bp*\bp))/2}
    \pgfmathsetmacro{\eb}{sqrt(\ea*\ea - pi*pi/4)}
    \draw[col_lgf2, dashed] ({pi/2},0) ellipse ({\ea} and {\eb});
    \node[col_lgf2!80!black, above] at ({pi/2},{\eb}) {Bernstein ellipse through $\pm\ii\eta$};
    \draw[->, black!60] (-0.55,0) -- (3.9,0) node[right] {$\operatorname{Re} k$};
    \draw[->, black!60] (0,{-\cut-0.15}) -- (0,{\cut+0.3}) node[right] {$\operatorname{Im} k$};
    \draw[very thick] (0,0) -- (3.1416,0);
    \node[below] at (3.1416,-0.42) {$\pi$};
    \foreach \s in {1,-1} {
        \draw[col_lgf8, very thick] (0,{\s*\bp}) -- (0,{\s*\cut});
        \node[col_lgf8] at (0,{\s*\bp}) {$\times$};
    }
    \foreach \x in {0.0290,0.1506,0.3614,0.6482,0.9930,1.3741,1.7675,2.1486,2.4934,2.7802,2.9910,3.1126} {\fill (\x,0.09) circle (1.1pt);}
    \foreach \x in {0.0000,0.2618,0.5236,0.7854,1.0472,1.3090,1.5708,1.8326,2.0944,2.3562,2.6180,2.8798,3.1416} {\draw (\x,-0.09) circle (1.1pt);}
    \node[right, black!70] at (3.6,0.62) {Gauss--Legendre};
    \node[right, black!70] at (3.6,-0.62) {trapezoidal};
    \draw[black!50] (3.55,0.55) -- (3.2,0.14);
    \draw[black!50] (3.55,-0.55) -- (3.25,-0.14);
    \node[col_lgf8, left] at (-0.08,0.5) {$\ii\eta$};
    \node[col_lgf8, left] at (-0.08,-0.5) {$-\ii\eta$};
    \draw[->, col_lgf8!80!black] (-0.4,\bp) -- (-0.4,0.04);
    \node[left, col_lgf8!80!black] at (-0.45,{\bp/2}) {$k_z\to0$};
\end{tikzpicture}
    \caption{Singularities of $I(n_x;k,k_z)$ in the complex $k_y$ plane at $k_z = \pi/8$, for any of the stencils considered: branch points $\pm\ii\eta$, $\eta\simeq k_z$, with their cuts (truncated), strip of analyticity (shaded), Bernstein ellipse through $\pm\ii\eta$ (dashed), and the interval of integration (thick) with the nodes of the $12$-point Gauss--Legendre (dots) and trapezoidal (circles) rules.}
    \label{fig:singularities}
\end{figure}

\paragraph{Step 2: the $k_y$ integral by Gauss--Legendre quadrature on $[0,\pi]$}
The integrand of the remaining integral is even in $k_y$, so
\begin{equation}
    G_{k_z}(n_x,n_y) = \frac{h^2}{\pi}\int_0^{\pi}\cos(n_y k)\,I(n_x;k,k_z)\dd k
    \approx \frac{h^2}{\pi}\sum_{q=1}^{Q} w_q\cos(n_y k_q)\,I(n_x;k_q,k_z).
    \label{eq:outer}
\end{equation}
For $k_z\neq0$, $e_0^2-e_1^2>0$ for every real $k$, so $I(n_x;k,k_z)$ and the integrand of \eqqref{eq:outer} are analytic in a strip about the real axis, and the quadrature error decays exponentially with the number of nodes $Q$. The rate is set by the nearest singularity, the branch point of \eqqref{eq:residue} at $k = \pm\ii\eta$, $\eta = 2\,\mathrm{arcsinh}\sin(k_z/2)\simeq k_z$, which closes onto the endpoint $k = 0$ as $k_z\to0$ (\fref{fig:singularities}).

How much that costs depends on the rule. The trapezoidal rule on $[0,\pi]$ converges exponentially~\cite{trefethen_exponentially_2014}, as established for the seven-point integrand in~\cite[Thm.~5.2]{hou_fast_2025}, but its rate is set by the width of the strip: the error decays like $\eul^{-2|k_z|Q}$ and machine precision needs $Q\simeq16/|k_z|$ nodes. The Gauss--Legendre rate is set instead by the Bernstein ellipse of foci $0$ and $\pi$ through the singularity~\cite[Thm.~19.3]{trefethen_approximation_2013}. That ellipse stays wide at mid-interval however close the branch point sits to the endpoint, where the nodes cluster. The exponent thus scales as $\sqrt{|k_z|}\,Q$, and machine precision needs $Q\simeq20/\sqrt{|k_z|}$ nodes (\tref{tab:quadrature}). For the seven-point stencil at $k_z = 2\pi/128$, computing $G_{k_z}(3,5)$ to machine precision takes $96$ Gauss--Legendre nodes against $320$ trapezoidal ones; at $k_z = \pi/2$, where the strip is wide, both need at most $32$.

The factor $\cos(n_y k)$ adds a second constraint: it oscillates $N_y/2$ times on $[0,\pi]$, which requires $Q\gtrsim 2N_y$. Each mode therefore uses
\begin{equation}
    Q(k_z) = \max\left(2N_y,\ \lceil 40/\sqrt{|k_z|}\rceil\right) + 32,
    \label{eq:Q}
\end{equation}
as a single Gauss--Legendre rule on $[0,\pi]$, i.e.\ twice the constant found above plus $32$ nodes, as a margin. The second term only exceeds the first for $|k_z| < (20/N_y)^2$, i.e.\ on the lowest modes of long periodic directions. The nodes and weights are generated to a few ulps at any order (Appendix~\ref{app:numerics}), and machine precision is observed in every case (\sref{sec:valid_kernel}).

The slice $k_z = 0$ would also be amenable to \eqqref{eq:outer} after subtraction of the integrand at $\bn = 0$, which restores analyticity and fixes the additive constant by $G_0(0,0) = 0$; this is the route taken for the second-order stencil in~\cite{hou_adaptive_2024}, where the subtracted value is itself available in closed form. We keep the tabulated two-dimensional LGF there, which is exact and grid independent.

\paragraph{Step 3: Hockney--Eastwood assembly}
For each mode $k_z\neq0$, $G_{k_z}$ on $[0,N_x]\times[0,N_y]$ is mirrored into a $2N_x\times2N_y$ even-periodic array and transformed; the result is real and replaces the $k_z\neq 0$ case of \eqqref{eq:formula}. Since $G_{k_z}$ is even in $k_z$, only $N_z/2+1$ slices are distinct. The construction applies verbatim to a symmetric (even or odd) instead of periodic third direction, whose modes are $k_z = (j_z+\varsigma)\pi/N_z$ with $\varsigma\in\{0,\tfrac12\}$: $k_z$ enters \eqqref{eq:outer} as a real parameter only.

\subsection{Regularised kernels}
\label{sec:regularised}
The construction is not tied to a lattice operator. The regularised kernels of~\cite{winckelmans1993contributions, hejlesen_high_2013} (\hej{2}--\hej{10}) are the free-space Green's function of the continuous Laplacian convolved with a mollifier of width $\sigma = \alpha h$ ($\alpha = 2$ throughout). Their extension to a periodic direction by~\cite{juul_spietz_regularization_2018} is \eqqref{eq:formula} with the regularised symbol in the $k_z\neq 0$ case. \Sref{sec:error} therefore applies to them, with one addition. Their symbol $\Ghat_\sigma$ is the transform of a continuum kernel and, unlike a lattice symbol, is not $2\pi$-periodic in $\bk$. The kernel $G_\sigma$ restricted to the grid points has, by Poisson summation, the periodic symbol $\sum_{\boldsymbol{\ell}}\Ghat_\sigma(\bk+2\pi\boldsymbol{\ell})$; the Spietz kernel samples $\Ghat_\sigma$ itself at the grid frequencies, i.e.\ keeps the term $\boldsymbol{\ell} = 0$ alone. The resulting error decays with $k_z$ like the transform of the mollifier, about $\exp(-(\alpha k_z)^2/2)$, not through screening.

For $k_z\neq0$ the regularised slice kernel is one Fourier integral away from the regularised one-unbounded/two-periodic kernel $G_1(x;k,\sigma)$ of~\cite{juul_spietz_regularization_2018}. $G_1$ is known in closed form and, the mollifier being isotropic, depends on the two spectral wavenumbers only through their magnitude (positions and wavenumbers are dimensional here, the kernel being a continuum one):
\begin{equation}
    G_{k_z}(x,y) = \frac{1}{\pi}\int_0^{\infty}\cos(k_y y)\,G_1\!\left(x;\sqrt{k_y^2+k_z^2},\sigma\right)\dd k_y .
    \label{eq:outer_hej}
\end{equation}
The integrand decays like $\eul^{-(k\sigma)^2/2}$, so the range is truncated at $k_{\max} = 9/\sigma$, below machine precision, and \eqqref{eq:outer_hej} is evaluated by the same Gauss--Legendre rule as \eqqref{eq:outer}, on $[0,k_{\max}]$, with
\begin{equation}
    Q(k_z) = \max\left(\left\lceil \frac{2k_{\max}L_y}{\pi}\right\rceil,\ \left\lceil 40\sqrt{\frac{k_{\max}}{\pi|k_z|}}\,\right\rceil\right) + \lceil 4k_{\max}\sigma\rceil + 64
    \label{eq:Q_hej}
\end{equation}
nodes, i.e.\ \eqqref{eq:Q} in the variable $\pi k/k_{\max}$ and with an extra term to resolve the Gaussian scale of the integrand. The two unbounded directions need not share a spacing. The mode $k_z = 0$ keeps its two-dimensional regularised kernel, as before.

\subsection{Implementation}
The kernels above are implemented in \flups{}~\cite{caprace_flups_2021, balty_flups_2023}, an open-source Fourier-based library of unbounded Poisson solvers for CPU and GPU clusters~\cite{flups}. The library carries its own distributed transforms, a component otherwise borrowed from libraries such as~\cite{pekurovsky_p3dfft_2012, gholami_accfft_2016, dalcin_fast_2019, krzhizhanovskaya_heffte_2020}. On top of them it provides the Hockney--Eastwood machinery and the lattice, regularised and continuous kernels of \sref{sec:lgf}, for any combination of unbounded, periodic and symmetric directions. It also contains the approximation \eqqref{eq:formula} that the present construction replaces. All the results of \sref{section:validation} are obtained with it. The slice assembly \eqqref{eq:outer} is carried out by blocks of modes, each block as one matrix product of the vendor BLAS (\sref{sec:valid_cost}). The lattice kernels require an isotropic lattice, since the closed forms of~\cite{gabbard_lattice_2023} are written for a single spacing; the regularised ones do not.


\section{Verification}
\label{section:validation}
All results are obtained with \flups{} on cell-centred grids of spacing $h$, with the two unbounded directions doubled by the HE algorithm and the third direction periodic unless stated otherwise, on the LUMI-G supercomputer (AMD MI250X, one MPI rank per graphics compute die, GCD). Errors are measured in the maximum norm $E_\infty = \Linf{u^h - u_{\mathrm{ref}}}$ over the original box. The kernel of \eqqref{eq:formula} is referred to as the \emph{Spietz kernel}, that of \sref{sec:exact} as \emph{exact}.

\subsection{Kernel-level verification}
\label{sec:valid_kernel}
The construction of \sref{sec:exact} is first checked on the slices themselves, independently of the solver, for the eight lattice kernels (\lgf{2}--\lgf{8}, \mehr{4L}, \mehr{6L}, \mehr{4F}, \mehr{6F}). Reference values of $G_{k_z}(n_x,n_y)$ are computed by two independent routes. The first is a brute-force inverse FFT of $\sigR/\sigL$ on a grid padded $32$ times, so that the image error \eqqref{eq:bound} is below machine precision. The second combines a dense trapezoidal rule in $k_x$ with the Gauss--Legendre rule of \eqqref{eq:outer} in $k_y$. For $N = 16$ and $k_z\in\{\pi/16,\pi/4,\pi/2,3\pi/4,\pi\}$, the two routes agree to $2\times10^{-15}$; the last two slices lie where the coefficient $e_1$ of \eqqref{eq:affine} is negative for \mehr{4} (Appendix~\ref{app:numerics}). \Tref{tab:lgdelta} reports two checks of the slices computed by \flups{} on the accelerator: their deviation from the references, and the residual of $\Lop G = \Rop\delta$ on each slice. Both are at machine precision for every kernel and every slice. The regularised kernels satisfy no lattice identity and have no independent tabulation; they are verified by the two solver-level tests below.

\begin{table}[htbp]
    \centering
    \caption{Kernel-level verification of the exact screened slices for $N = 16$, maximum over the five slices $k_z\in\{\pi/16,\pi/4,\pi/2,3\pi/4,\pi\}$: residual of $\Lop G = \Rop\delta$ with the true stencils and deviation from the reference values obtained by two independent quadratures.}
    \label{tab:lgdelta}
    \footnotesize
\begin{tabular}{l c c}
    \toprule
    kernel & $\max|\Lop G - \Rop\delta|$ & $\max|G - G_{\mathrm{ref}}|$ \\
    \midrule
    \lgf{2}   & $1\cdot10^{-15}$ & $1\cdot10^{-16}$ \\
    \lgf{4}   & $2\cdot10^{-15}$ & $5\cdot10^{-16}$ \\
    \lgf{6}   & $2\cdot10^{-15}$ & $1\cdot10^{-15}$ \\
    \lgf{8}   & $2\cdot10^{-15}$ & $3\cdot10^{-16}$ \\
    \mehr{4L} & $1\cdot10^{-15}$ & $6\cdot10^{-16}$ \\
    \mehr{6L} & $1\cdot10^{-15}$ & $3\cdot10^{-16}$ \\
    \mehr{4F} & $1\cdot10^{-15}$ & $4\cdot10^{-16}$ \\
    \mehr{6F} & $9\cdot10^{-16}$ & $3\cdot10^{-16}$ \\
    \bottomrule
\end{tabular}

\end{table}

\subsection{A padding test}
\label{sec:valid_pad}
\begin{figure}[htbp]
    \centering
\begin{tikzpicture}[x=0.95cm, y=0.95cm, >=stealth, font=\footnotesize]
    \def\m{0.7}
    \newcommand{\noisecell}[3]{
        \pgfmathsetseed{4}
        \foreach \i in {0,...,7} {
            \foreach \j in {0,...,7} {
                \pgfmathsetmacro{\v}{rand}
                \pgfmathsetmacro{\a}{abs(\v)*70}
                \ifdim\v pt>0pt
                    \fill[col_lgf2!\a, opacity=#3] ({#1+\i/8},{#2+\j/8}) rectangle ++(0.125,0.125);
                \else
                    \fill[col_lgf8!\a, opacity=#3] ({#1+\i/8},{#2+\j/8}) rectangle ++(0.125,0.125);
                \fi
            }
        }
    }
    \newcommand{\padpanel}[5]{
        \begin{scope}[shift={#1}]
            \pgfmathtruncatemacro{\P}{2*#2}
            \begin{scope}
                \clip (-\m,-\m) rectangle ({\P+\m},{\P+\m});
                \foreach \p in {-1,0,1} {
                    \foreach \q in {-1,0,1} {
                        \ifnum\numexpr\p*\p+\q*\q\relax>0
                            \noisecell{\p*\P}{\q*\P}{0.45}
                            \draw[black!40, dashed] ({\p*\P},{\q*\P}) rectangle ++(\P,\P);
                        \fi
                    }
                }
            \end{scope}
            \fill[white] (0,0) rectangle (\P,\P);
            \fill[pattern=north east lines, pattern color=black!25] (0,0) rectangle (\P,\P);
            \fill[white] (0,0) rectangle (#2,#2);
            \ifnum#2>1
                \fill[col_lgf2!5] (0,0) rectangle (#2,#2);
                \draw (0,0) rectangle (#2,#2);
            \fi
            \noisecell{0}{0}{1}
            \draw[very thick] (0,0) rectangle (1,1);
            \draw (0,0) rectangle (\P,\P);
            \node[above, align=center] at ({\P/2},{\P+\m}) {#3\\[-1pt]{\scriptsize\color{black!60}Spietz kernel: images at period $#4$}};
            \draw[<->] (1,0.5) -- (\P,0.5) node[midway, above=1pt, fill=white, inner sep=1pt, font=\scriptsize] {$d\ge#5$};
        \end{scope}
    }
    \padpanel{(0,2)}{1}{(a) box, $N\times N$}{2L_x}{L_x}
    \padpanel{(5.2,0)}{2}{(b) enlarged box, $2N\times2N$, same $h$}{4L_x}{3L_x}
    \node[align=center, font=\scriptsize] at (1.5,3.5) {zero\\padding};
    \draw[<-] (-0.05,2.5) -- (-0.85,2.5) node[left, align=right, font=\scriptsize, inner sep=1pt] {white\\noise};
    \draw[|<->|] (0,1.88) -- (1,1.88) node[pos=0.3, below, font=\scriptsize] {$L_x$};
    \node[align=center, font=\scriptsize] at (5.2+1.5,1.45) {no\\source};
    \node[align=center, font=\scriptsize] at (5.2+3,3) {zero padding};
    \draw[->] (9.25,-1.35) -- (9.65,-1.35) node[right] {$x$};
    \draw[->] (9.25,-1.35) -- (9.25,-0.95) node[above] {$y$};
    \node[draw, rounded corners, fill=white, align=center, inner sep=3pt, font=\scriptsize] (cmp) at (0.95,-0.2) {
        on the common box (thick), mode $j$ along $z$:\\[2pt]
        $\Delta\hat u_j = \hat u^{(a)}_j-\hat u^{(b)}_j$\\[1pt]
        $e_j = \max|\Delta\hat u_j|\,/\max|\hat u|$\\[3pt]
        exact kernel: both free-space, $e_j\approx\epsilon$\\
        Spietz kernel: the image fields differ
    };
    \draw[->, thick] (0.8,1.95) -- (cmp.north -| 0.8,0);
    \draw[->, thick] (5.2+0.5,-0.02) to[out=-100, in=0] (cmp.east);
\end{tikzpicture}
    \caption{Padding test in the $(x,y)$ plane: the same white-noise source solved on the box (a) and on the box enlarged twice with the same $h$ (b), compared on the common box.}
    \label{fig:padtest}
\end{figure}
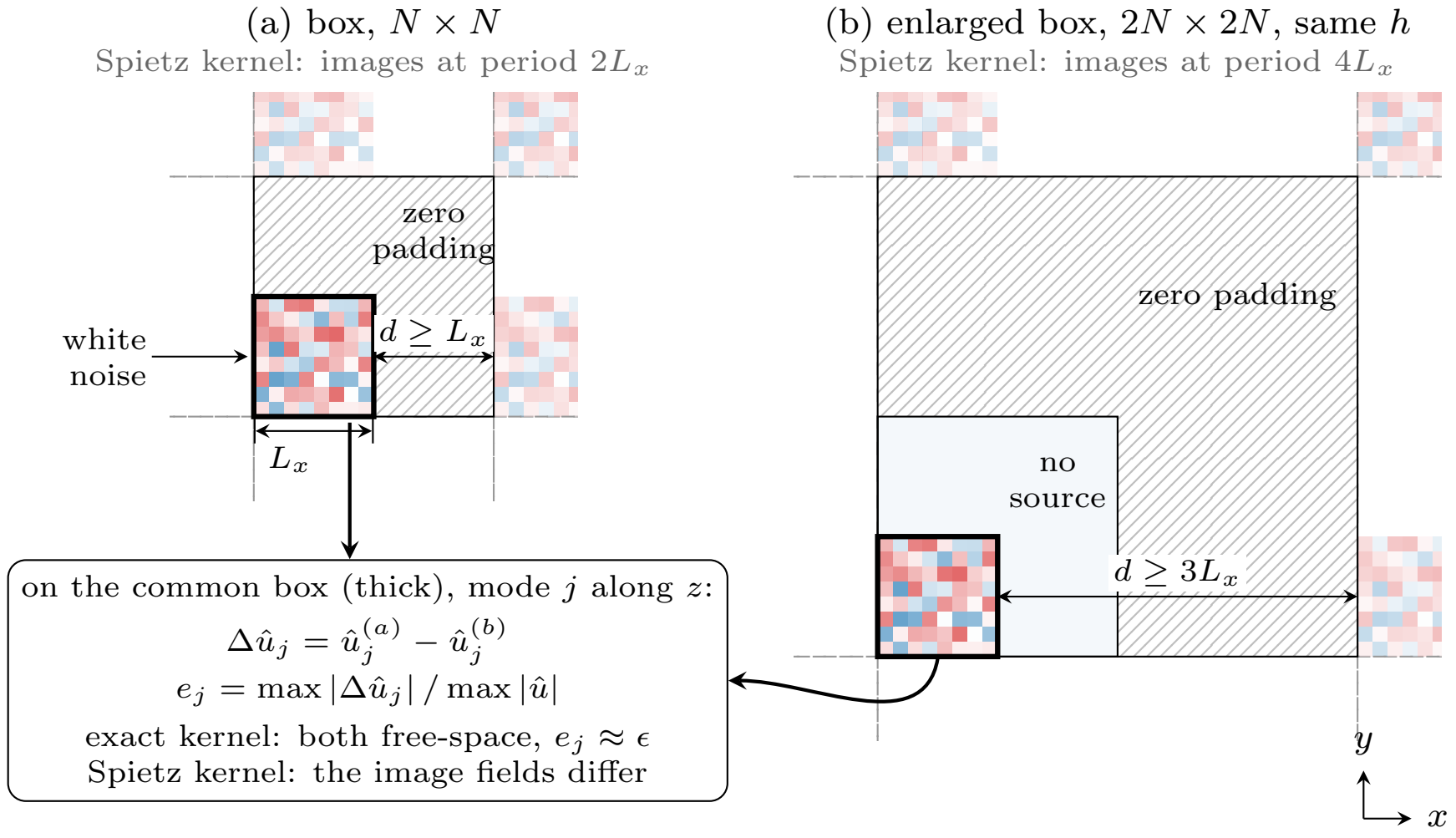

Following \sref{sec:error}, the first solver-level test uses a compact \emph{source}: white noise, uniform in $[-1,1]$ at every point of the box, so that every mode of the periodic direction, hence every slice of the kernel, is exercised. The same source is solved, with the same $h$, on the box and on an enlarged box twice as large in the two unbounded directions, and the two solutions are compared on the common box, mode by mode along $z$ (\fref{fig:padtest}). For an exact free-space kernel they coincide; for the Spietz kernel they differ by the field of the images, \eqqref{eq:split}.

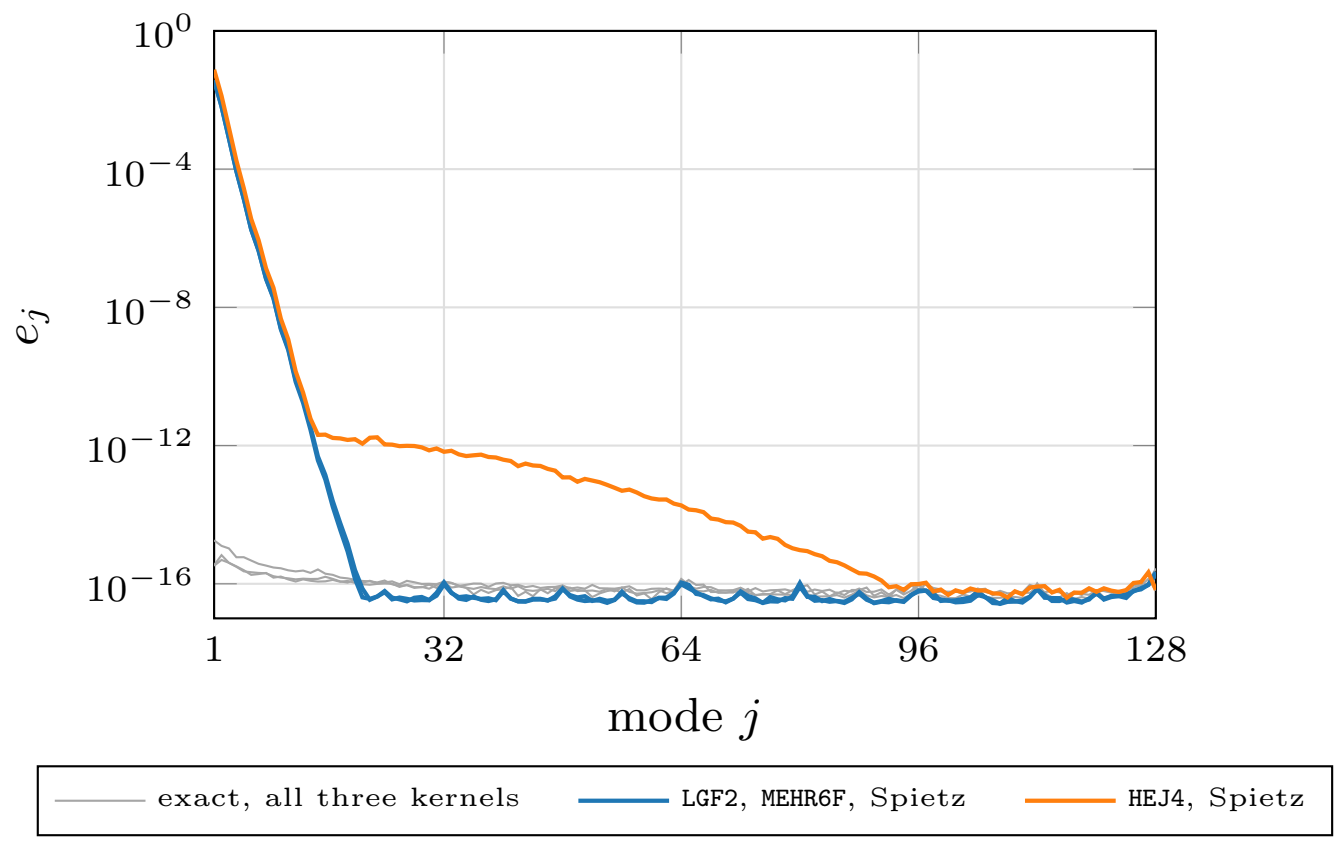
\begin{figure}[htbp]
    \centering
\begin{tikzpicture}
    \begin{axis}[
        lgfplot,
        width=0.6\textwidth, height=0.42\textwidth,
        ymode=log,
        xlabel={mode $j$}, ylabel={$e_j$},
        xmin=1, xmax=128, xtick={1,32,64,96,128},
        restrict x to domain=1:128,
        ymin=1e-17, ymax=1e0,
        ytick={1e0,1e-4,1e-8,1e-12,1e-16},
        legend style={anchor=north, at={(0.5,-0.25)}, legend columns=3, font=\tiny, /tikz/every even column/.append style={column sep=0.3cm}},
        legend cell align=left,
        table/header=false,
    ]
        \addplot[gray!70, thin]                table[x index=0, y index=5] {data/padnoise/per_qs_N64_Lz4_LGF_2.dat};
        \addlegendentry{exact, all three kernels}
        \addplot[gray!70, thin, forget plot]   table[x index=0, y index=5] {data/padnoise/per_qs_N64_Lz4_MEHR_6F.dat};
        \addplot[gray!70, thin, forget plot]   table[x index=0, y index=5] {data/padnoise/per_qs_N64_Lz4_HEJ_4.dat};
        \addplot[col_lgf2]                     table[x index=0, y index=5] {data/padnoise/per_samp_N64_Lz4_LGF_2.dat};
        \addlegendentry{{\lgf{2}, \mehr{6F}, Spietz}}
        \addplot[col_lgf2, forget plot]        table[x index=0, y index=5] {data/padnoise/per_samp_N64_Lz4_MEHR_6F.dat};
        \addplot[col_lgf4]                     table[x index=0, y index=5] {data/padnoise/per_samp_N64_Lz4_HEJ_4.dat};
        \addlegendentry{{\hej{4}, Spietz}}
    \end{axis}
\end{tikzpicture}
    \caption{Padding test with a white-noise source, $N = 64$, $L_z = 4L_x$, LUMI-G: error $e_j = \max|\Delta\hat u_j|/\max|\hat u|$ of the mode $j$ ($k_z = 2\pi j/N_z$), i.e.\ the largest change of that mode of the solution when the box is enlarged, relative to its largest mode, for \lgf{2}, \mehr{6F} and \hej{4}, Spietz kernel (\lgf{2} and \mehr{6F} coincide) and exact kernels (grey). The mode $j = 0$, exact for every kernel, is not shown.}
    \label{fig:padnoise}
\end{figure}

\Tref{tab:padcheck} reports the largest error over the modes, $\max_j e_j$, relative to the solution as in \fref{fig:padnoise}. The exact kernels are at machine precision, about $10^{-15}$. The Spietz kernel grows with the aspect ratio, from $10^{-7}$ to $0.15$ for the lattice kernels, in line with the estimates of \sref{sec:error} (only \lgf{2} stands apart on the flat box, by the dispersion of its stencil). The regularised kernels are about three times higher only because they smooth the noise, so their solution is smaller: the difference itself is the same for every kernel. \Fref{fig:padnoise} shows that this error sits in the lowest modes and decays exponentially with $j$; the mode $k_z = 0$ is exact. The regularised kernels are affected on many more modes, up to $j\simeq80$--$100$ against $20$ for the lattice ones, by the truncation of their symbol (see \sref{sec:regularised}). The same holds on node-centred grids and with a symmetric (even/odd) third direction, where the Spietz kernel is worse, its lowest mode $k_z = \pi/N_z$ being screened half as much.

\begin{table}[htbp]
    \centering
    \caption{Padding test with a white-noise source on $N_x = N_y = 16$ (and $N_z = 16\,L_z/L_x$): largest error over the modes, $\max_j e_j$, with $e_j$ as in \fref{fig:padnoise}, Spietz versus exact kernels, lattice (\lgf{2}--\lgf{8}, \mehr{4L/6L/4F/6F}) and regularised (\hej{2}--\hej{10}), for aspect ratios $L_z/L_x$ from $1/2$ to $8$.}
    \label{tab:padcheck}
    \scriptsize
\setlength{\tabcolsep}{4pt}
\begin{tabular}{l l c c c c c}
    \toprule
    & $L_z/L_x$ & $1/2$ & $1$ & $2$ & $4$ & $8$ \\
    \midrule
    Spietz kernel & \lgf{2} & $1.3\cdot10^{-7}$ & $3.4\cdot10^{-5}$ & $1.5\cdot10^{-3}$ & $5.3\cdot10^{-2}$ & $1.5\cdot10^{-1}$ \\
     & other lattice (max) & $7.2\cdot10^{-8}$ & $3.2\cdot10^{-5}$ & $1.5\cdot10^{-3}$ & $5.4\cdot10^{-2}$ & $1.5\cdot10^{-1}$ \\
     & regularised (max) & $3.4\cdot10^{-7}$ & $9.0\cdot10^{-5}$ & $4.8\cdot10^{-3}$ & $1.7\cdot10^{-1}$ & $4.3\cdot10^{-1}$ \\
    \midrule
    exact slices & lattice (max) & $1.3\cdot10^{-15}$ & $1.1\cdot10^{-15}$ & $1.1\cdot10^{-15}$ & $1.3\cdot10^{-15}$ & $1.1\cdot10^{-15}$ \\
     & regularised (max) & $1.3\cdot10^{-15}$ & $5.1\cdot10^{-16}$ & $6.5\cdot10^{-16}$ & $9.3\cdot10^{-16}$ & $1.1\cdot10^{-15}$ \\
    \bottomrule
\end{tabular}

\end{table}

\subsection{An axially modulated tube source}
\label{sec:valid_wavy}
The second solver-level test has an analytical solution that is $\mathcal{O}(1)$ on the unbounded faces and an exactly compact source, so that the convergence at the order of each stencil can be observed. Its radial profile is borrowed from the vortex tube of~\cite{winckelmans2004encyclopedia}, a tube of radius $R$ aligned with the periodic direction, whose vorticity is the compactly supported, $C^\infty$ bump
\begin{equation}
    \omega(r) = \frac{1}{\pi R^2 E_2(1)}\exp\left(\frac{-1}{1-(r/R)^2}\right),\quad r<R,\qquad \omega = 0 \text{ otherwise},
\end{equation}
of unit circulation ($E_2$ is the exponential integral). A straight tube has no content on the modes $k_z\neq0$ and cannot exercise the kernels of interest; we therefore modulate it along its axis and solve, in the scalar setting of \eqqref{eq:poisson},
\begin{equation}
    f(x,y,z) = \omega(r)\sum_{j=1}^{2} a_j\cos(k_j z - \varphi_j),\qquad k_j = \frac{2\pi j}{L_z},\quad a = (1,\tfrac12),\quad \varphi = (0,\tfrac{\pi}{2}),
    \label{eq:wavy_rhs}
\end{equation}
with $r$ the distance to the axis at the centre of the transverse section and $R = 0.2L_x$. The mode $k_z = 0$ is left out on purpose: it is the classical two-dimensional problem, whose kernel is unchanged, and it carries the additive constant of the two-dimensional LGF. \Eqqref{eq:wavy_rhs} is a manufactured scalar source, not a vorticity field; it has what the test needs: an exactly compact source, a solution that is $\mathcal{O}(1)$ on the unbounded faces, and controlled energy on the modes $k_z\neq0$.

The exact solution separates, $u = \sum_j a_j\cos(k_j z-\varphi_j)\,\psi_j(r)$, where each $\psi_j$ solves the screened radial problem $-(\nabla^2_{\mathrm{2D}} - k_j^2)\psi_j = \omega$, i.e.\ with the radial Green's function $I_0(k r_<)K_0(k r_>)$ of a source ring at radius $s$~\cite[\S3.11]{jackson_classical_1999}, $r_< = \min(r,s)$ and $r_> = \max(r,s)$,
\begin{equation}
    \psi_j(r) = K_0(k_j r)\int_0^r I_0(k_j s)\,\omega(s)\,s\dd s + I_0(k_j r)\int_r^R K_0(k_j s)\,\omega(s)\,s\dd s,
    \label{eq:wavy_sol}
\end{equation}
so that $\psi_j = A_j K_0(k_j r)$ outside the core: the solution decays like $\eul^{-k_j r}$ and is $\mathcal{O}(1)$ on the faces of any box a few radii wide, while the source vanishes identically there. The reference solution is therefore semi-analytic; its radial integrals are evaluated to $4\times10^{-16}$ (Appendix~\ref{app:numerics}), a factor of three below the smallest error of \fref{fig:wavy_convergence}, so that the curves measure the discretisation error of the stencils rather than the accuracy of the reference.

On a cubic box the lowest mode is screened by $\eul^{-2\pi\cdot 1.5}\approx10^{-4}$ before it reaches the images and the image error of \eqqref{eq:bound} would only be visible at high order; the box is therefore taken four times longer than wide, $L_z = 4L_x$, where $k_1L_x = \pi/2$ and the image error is $\mathcal{O}(10^{-1})$.

\Fref{fig:wavy_convergence} shows $E_\infty$ against the transverse resolution $N = N_x = N_y$ ($N_z = 4N$) for the six lattice kernels (left) and the five regularised ones (right), from $N = 32$ to $N = 4096$, i.e.\ up to $4096^2\times16384\simeq2.7\times10^{11}$ unknowns on $4096$ GCDs. With the Spietz kernel the error does not converge at all: it stays at $3\times10^{-2}$, the same for every kernel, and is drawn once per panel. Every exact lattice kernel converges at the order of its stencil until it meets the round-off floor of the solve, about $10^{-15}$, at the level reported in~\cite{gabbard_lattice_2023} for the one-unbounded kernels. The regularised kernels converge as well, at orders that approach their nominal values from below as the mollifier width $\sigma = 2h$ shrinks, down to the same round-off floor.

\begin{figure}[htbp]
    \centering
\begin{tikzpicture}
    \begin{groupplot}[
        group style={group size=2 by 1, horizontal sep=0.4cm, ylabels at=edge left, yticklabels at=edge left},
        lgfplot, width=0.52\textwidth, height=0.42\textwidth,
        xmode=log, ymode=log,
        log basis x=2,
        xlabel={$L_z/L_x$}, ylabel={$E_\infty$},
        xtick={1,2,4,8,16,32,64},
        xticklabels={1,2,4,8,16,32,64},
        ymin=1e-8, ymax=1e2,
        legend cell align=left,
    ]
    \nextgroupplot[title={\small single tube}, legend to name=leg:wavy_aspect,
        legend style={legend columns=3, font=\tiny, /tikz/every even column/.append style={column sep=0.3cm}}]
        \addplot[black!18, line width=4pt, mark=none, forget plot] table[x index=0, y index=1] {data/wavy_aspect_image_model.dat}
            node[pos=0.8, below right, font=\tiny, black!60] {image field};
        \addplot[col_lgf2,  mark=o,   dashed]      table[x index=0, y index=2] {data/wavy_aspect_LGF_2_old_lumi.dat};
        \addplot[col_mehr6, mark=square, dashed]   table[x index=0, y index=2] {data/wavy_aspect_MEHR_6F_old_lumi.dat};
        \addplot[col_lgf4,  mark=triangle, dashed] table[x index=0, y index=2] {data/wavy_aspect_HEJ_4_old_lumi.dat};
        \addplot[col_lgf2,  mark=*]                table[x index=0, y index=2] {data/wavy_aspect_LGF_2_new_lumi.dat};
        \addplot[col_mehr6, mark=square*]          table[x index=0, y index=2] {data/wavy_aspect_MEHR_6F_new_lumi.dat};
        \addplot[col_lgf4,  mark=triangle*]        table[x index=0, y index=2] {data/wavy_aspect_HEJ_4_new_lumi.dat};
        \legend{{\lgf{2}, Spietz}, {\mehr{6F}, Spietz}, {\hej{4}, Spietz}, {\lgf{2}, exact}, {\mehr{6F}, exact}, {\hej{4}, exact}}
    \nextgroupplot[title={\small counter-rotating pair}]
        \addplot[black!18, line width=4pt, mark=none, forget plot] table[x index=0, y index=1] {data/wavy_pair_aspect_image_model.dat}
            node[pos=0.8, above, font=\tiny, black!60] {image field};
        \addplot[col_lgf2,  mark=o,   dashed]      table[x index=0, y index=2] {data/wavy_pair_aspect_LGF_2_old_lumi.dat};
        \addplot[col_mehr6, mark=square, dashed]   table[x index=0, y index=2] {data/wavy_pair_aspect_MEHR_6F_old_lumi.dat};
        \addplot[col_lgf4,  mark=triangle, dashed] table[x index=0, y index=2] {data/wavy_pair_aspect_HEJ_4_old_lumi.dat};
        \addplot[col_lgf2,  mark=*]                table[x index=0, y index=2] {data/wavy_pair_aspect_LGF_2_new_lumi.dat};
        \addplot[col_mehr6, mark=square*]          table[x index=0, y index=2] {data/wavy_pair_aspect_MEHR_6F_new_lumi.dat};
        \addplot[col_lgf4,  mark=triangle*]        table[x index=0, y index=2] {data/wavy_pair_aspect_HEJ_4_new_lumi.dat};
    \end{groupplot}
    \path (group c1r1.outer south) -- (group c2r1.outer south)
        node[pos=0.5, anchor=north, yshift=-2pt] {\pgfplotslegendfromname{leg:wavy_aspect}};
\end{tikzpicture}
    \caption{Axially modulated tube source (left) and counter-rotating pair of such tubes (right) at fixed transverse resolution $N = 128$, LUMI-G: maximum error against the aspect ratio $L_z/L_x$ for \lgf{2}, \mehr{6F} and \hej{4}, Spietz kernel (dashed) and exact kernels (solid); the grey bands are the field of the periodic images of the source on the doubled box, computed from the continuous solution (see text).}
    \label{fig:wavy_aspect}
\end{figure}
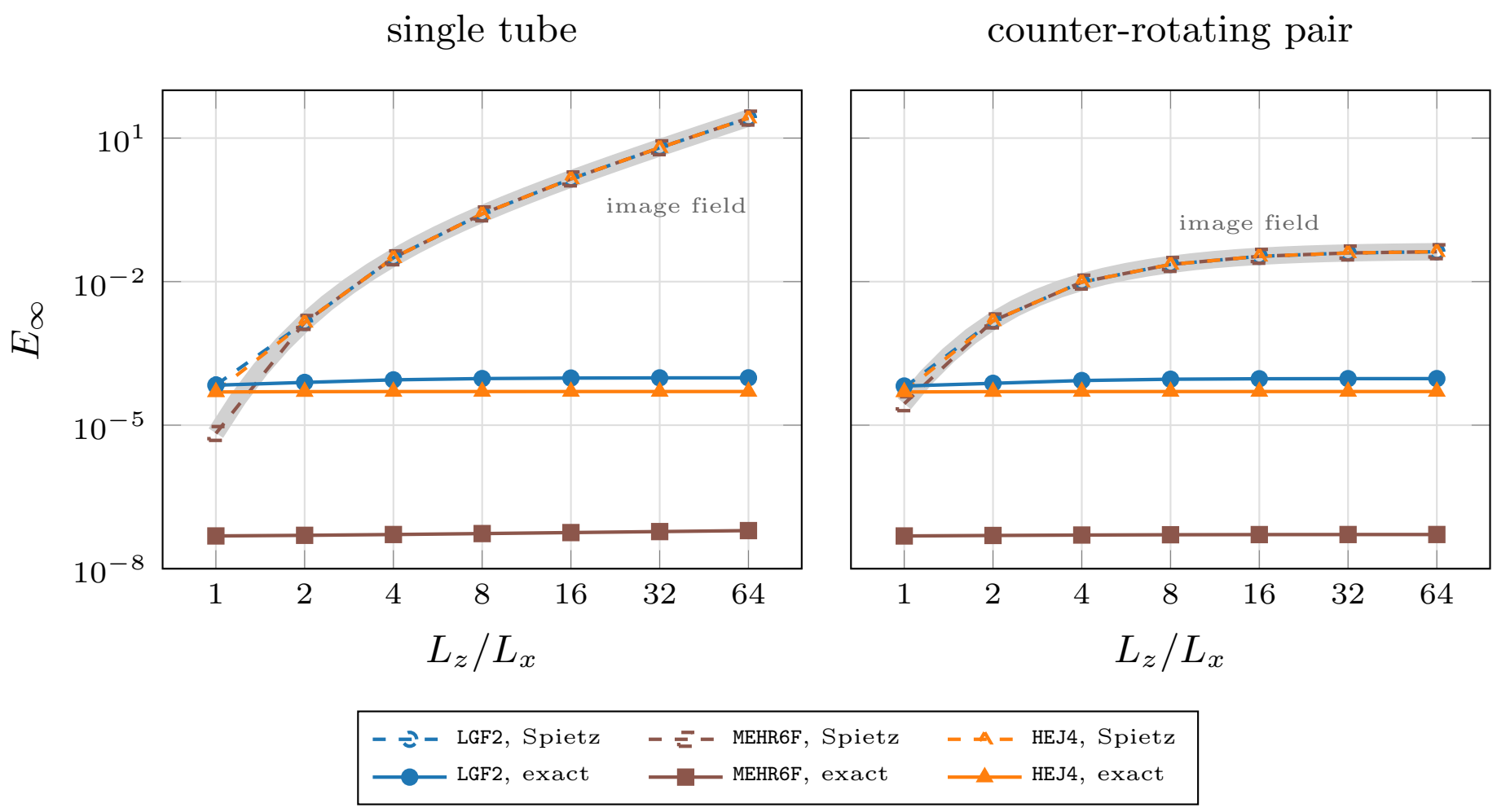

\Fref{fig:wavy_aspect} isolates the dependence on the aspect ratio at fixed $N = 128$, from $L_z/L_x = 1$ to $64$. The error of the exact kernels is the discretisation error and varies by at most a factor $1.4$ over that range. The error of the Spietz kernel is the same for the three kernels from $L_z/L_x = 2$ on, and reaches $26$ at $L_z/L_x = 64$, far more than the solution itself. On the cube alone it is still hidden under the discretisation error of \lgf{2} and \hej{4}, though already $140$ times that of \mehr{6F}. It is the field of the periodic images, \eqqref{eq:split}. Outside the core, each mode of the solution is $A_jK_0(k_jr)$; summed over the images of the doubled box, it gives the grey band. The band matches the data to $2\%$ wherever it exceeds the discretisation error, and to $0.03\%$ from $L_z/L_x = 16$ on. Both regimes of \sref{sec:error} appear: exponential growth well below $L_z = 4\pi L_x$ and, beyond, the algebraic growth of \eqqref{eq:bound_long}, since the mode $j = 1$ carries a unit net source.

The right panel repeats the sweep for a counter-rotating pair, the tube and its opposite at $x = L_x/4$ and $3L_x/4$, whose net source vanishes. The error of the Spietz kernel grows exponentially while the lowest mode is screened, as for the single tube, but levels off beyond instead of growing algebraically: it is $4.2\times10^{-2}$ at $L_z/L_x = 64$, against $26$ for the single tube, within $5\%$ of the plateau $|\boldsymbol{D}|/(16L_x)$ of \sref{sec:error}, and the image field of the pair matches it to $0.15\%$ or better.

\begin{figure}[htbp]
    \centering
\begin{tikzpicture}
    \begin{groupplot}[
        group style={group size=2 by 1, horizontal sep=0.4cm, ylabels at=edge left, yticklabels at=edge left},
        lgfplot, width=0.52\textwidth, height=0.42\textwidth,
        xmode=log, ymode=log,
        log basis x=2,
        xlabel={$N$}, ylabel={$E_\infty$},
        xtick={32,64,128,256,512,1024,2048,4096},
        xticklabels={32,64,128,256,512,1024,2048,4096},
        x tick label style={font=\tiny},
        ymin=1e-15, ymax=1e0,
        ytick={1e0,1e-2,1e-4,1e-6,1e-8,1e-10,1e-12,1e-14},
        xmin=32, xmax=4096, enlarge x limits=0.08,
        legend style={anchor=north, at={(0.5,-0.24)}, legend columns=3, font=\tiny, /tikz/every even column/.append style={column sep=0.15cm}},
        legend cell align=left,
    ]
    \nextgroupplot[title={\small lattice kernels}]
        \addplot[black!45, thick, mark=*, mark size=1.3pt, forget plot] table[x index=0, y index=2] {data/wavy_conv_LGF_2_old_lumi.dat};
        \addplot[col_lgf2,  mark=o]        table[x index=0, y index=2] {data/wavy_conv_LGF_2_new_lumi.dat};
        \addplot[col_lgf4,  mark=square]   table[x index=0, y index=2] {data/wavy_conv_LGF_4_new_lumi.dat};
        \addplot[col_lgf6,  mark=triangle] table[x index=0, y index=2] {data/wavy_conv_LGF_6_new_lumi.dat};
        \addplot[col_lgf8,  mark=diamond]  table[x index=0, y index=2] {data/wavy_conv_LGF_8_new_lumi.dat};
        \addplot[col_mehr4, mark=x]        table[x index=0, y index=2] {data/wavy_conv_MEHR_4F_new_lumi.dat};
        \addplot[col_mehr6, mark=+]        table[x index=0, y index=2] {data/wavy_conv_MEHR_6F_new_lumi.dat};
        \legend{\lgf{2}, \lgf{4}, \lgf{6}, \lgf{8}, \mehr{4F}, \mehr{6F}}
        \node[font=\tiny, black!70, anchor=north east, inner sep=1pt] at (axis cs:4096,3.1e-1) {Spietz, all kernels};
        \addplot[orderline, domain=32:4096, samples=2] {2e-3*(32/x)^2} node[pos=1, right, font=\tiny] {$2$};
        \addplot[orderline, domain=32:4096, samples=2] {2e-3*(32/x)^4} node[pos=0.97, above, font=\tiny] {$4$};
        \addplot[orderline, domain=32:4096, samples=2] {2e-3*(32/x)^6} node[pos=0.85, above, font=\tiny] {$6$};
        \addplot[orderline, domain=32:4096, samples=2] {2e-3*(32/x)^8} node[pos=0.7, left, font=\tiny] {$8$};
    \nextgroupplot[title={\small regularised kernels}]
        \addplot[black!45, thick, mark=*, mark size=1.3pt, forget plot] table[x index=0, y index=2] {data/wavy_conv_HEJ_4_old_lumi.dat};
        \addplot[col_lgf2,  mark=o]        table[x index=0, y index=2] {data/wavy_conv_HEJ_2_new_lumi.dat};
        \addplot[col_lgf4,  mark=square]   table[x index=0, y index=2] {data/wavy_conv_HEJ_4_new_lumi.dat};
        \addplot[col_lgf6,  mark=triangle] table[x index=0, y index=2] {data/wavy_conv_HEJ_6_new_lumi.dat};
        \addplot[col_lgf8,  mark=diamond]  table[x index=0, y index=2] {data/wavy_conv_HEJ_8_new_lumi.dat};
        \addplot[col_hej10, mark=pentagon] table[x index=0, y index=2] {data/wavy_conv_HEJ_10_new_lumi.dat};
        \legend{\hej{2}, \hej{4}, \hej{6}, \hej{8}, \hej{10}}
        \addplot[orderline, domain=32:4096, samples=2] {2e-3*(32/x)^2} node[pos=1, right, font=\tiny] {$2$};
        \addplot[orderline, domain=32:4096, samples=2] {2e-3*(32/x)^4} node[pos=0.97, below, font=\tiny] {$4$};
        \addplot[orderline, domain=32:4096, samples=2] {2e-3*(32/x)^6} node[pos=0.8, below, font=\tiny] {$6$};
        \addplot[orderline, domain=32:4096, samples=2] {2e-3*(32/x)^8} node[pos=0.7, left, font=\tiny] {$8$};
        \addplot[orderline, domain=32:4096, samples=2] {2e-3*(32/x)^10} node[pos=0.565, left, font=\tiny] {$10$};
    \end{groupplot}
\end{tikzpicture}
    \caption{Axially modulated tube source, $L_z = 4L_x$, LUMI-G: maximum error against the transverse resolution $N$ ($N_z = 4N$) for the exact kernels, lattice (left) and regularised (right); the grey line is the error of the Spietz kernel, the same for every kernel to three digits; dashed lines indicate the orders $2$ to $10$.}
    \label{fig:wavy_convergence}
\end{figure}
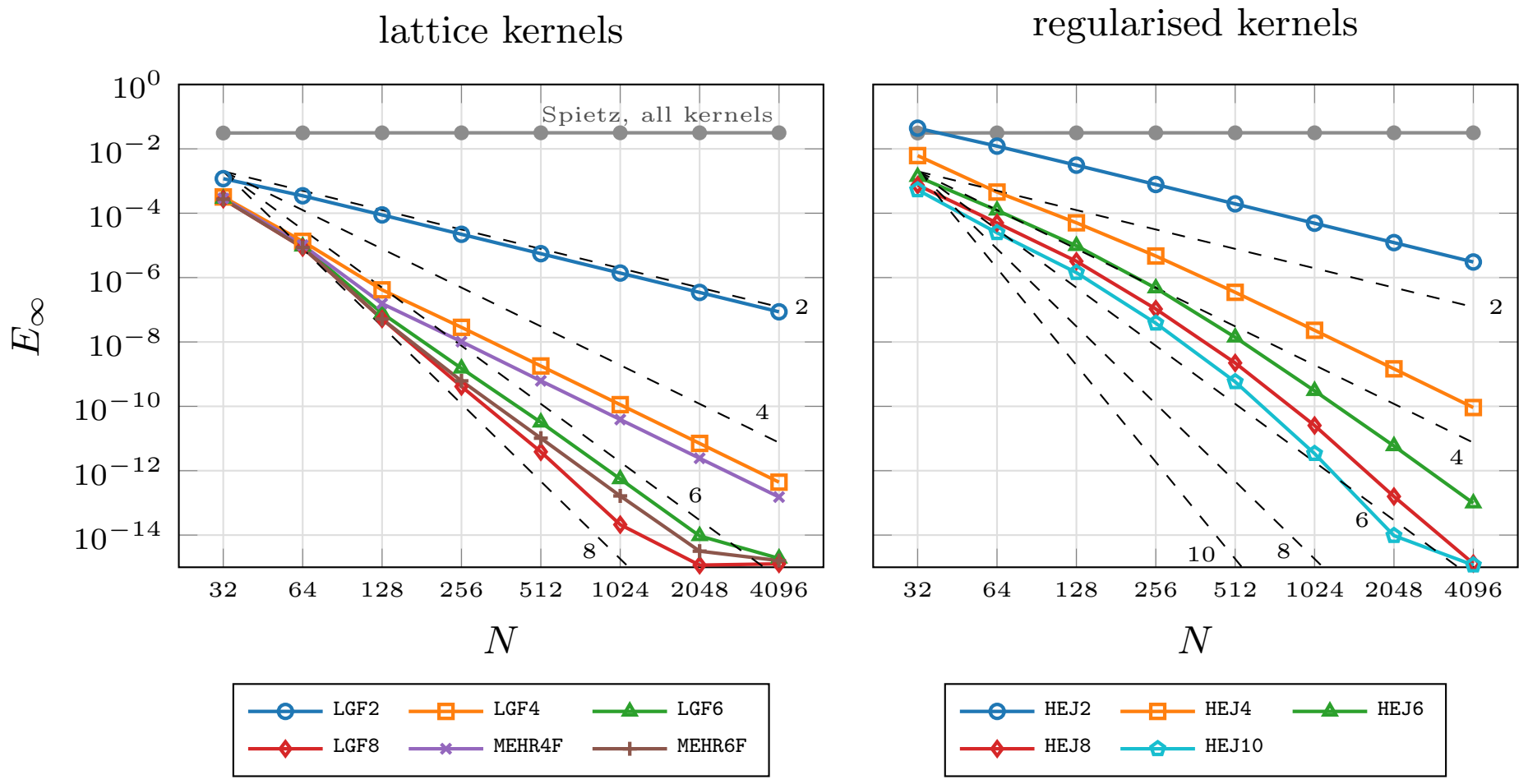

\subsection{Cost}
\label{sec:valid_cost}
\Tref{tab:cost} gives the cost of the exact kernels on LUMI-G, for the grids of \sref{sec:valid_wavy}, from $N = 512$ on $8$ GCDs to $N = 4096$ on $4096$ GCDs. Each kernel row is the time spent computing the kernel: the closed forms and the quadrature assembly. The transforms of the kernel and the creation of the FFT plans are listed apart. They do not depend on the kernel, and the Spietz kernel paid them as well, like every unbounded kernel of the library. The kernel rows are therefore the cost of exactness itself. The last row is the time of one solve on the same grid, given for reference.

The kernel is computed once per grid. It takes between $0.5$ and $0.9$~s on every grid of the sweep, that is between one fifth of a solve and one solve. At a fixed load per GCD (with the exception of $N = 1024$, where it is halved), the cost hardly grows with $N$: each rank assembles all the modes, but only on its own patch of the doubled plane. The whole setup, kernel included, is between one and six solves. The largest ratio occurs at $N = 1024$, where the halved load brings the solve down to half a second while the transforms and plans take $2.7$~s.

The multiply-adds of the assembly, \eqqref{eq:outer}, are $\mathcal{O}(N_zN^2Q)$ per grid but cost little. The slices are assembled by blocks of up to $64$ modes sharing the same rule: the values of the inner integral are stacked into one tall matrix. Then, a single product with the matrix of weighted cosines assembles the block. The vendor BLAS runs this shape close to the double-precision peak of a GCD, whereas one product per mode, a few tens of rows wide, runs far below it. What remains is dominated by the closed forms of Step~1, whose cost depends on the stencil: the eighth-order LGF is $40\%$ more expensive than the other kernels at $N = 4096$.

The regularised kernels use a $1.4$ times larger rule, \eqqref{eq:outer_hej} being cut off at $k_y\sigma = 9$, i.e.\ at $4.5/h$, against the half-period $\pi/h$ of the lattice symbol. The product absorbs it: both families cost the same to within $10\%$. The online solve, its memory footprint and its communication pattern are strictly unchanged.

\begin{table}[htbp]
    \centering
    \caption{Computation time of the Green's function, in seconds, for the grids of \sref{sec:valid_wavy} ($N_z = 4N$, doubled unbounded directions) on LUMI-G, one MPI rank per GCD, the rest of the Green's function setup (transforms of the kernel and creation of the FFT plans, the same for every kernel), and time of one solve on the same grid.}
    \label{tab:cost}
    \footnotesize
\begin{tabular}{l r r r r }
    \toprule
    $N$ & $512$ & $1024$ & $2048$ & $4096$ \\
    GCDs & $8$ & $128$ & $512$ & $4096$ \\
    \midrule
    \lgf{2} & 0.48 & 0.47 & 0.53 & 0.60 \\
    \lgf{4} & 0.50 & 0.49 & 0.56 & 0.67 \\
    \lgf{6} & 0.50 & 0.49 & 0.57 & 0.68 \\
    \lgf{8} & 0.54 & 0.56 & 0.70 & 0.94 \\
    \mehr{4F} & 0.50 & 0.47 & 0.55 & 0.64 \\
    \mehr{6F} & 0.49 & 0.48 & 0.55 & 0.65 \\
    \hej{2} & 0.48 & 0.47 & 0.54 & 0.62 \\
    \hej{4} & 0.48 & 0.47 & 0.54 & 0.63 \\
    \hej{6} & 0.49 & 0.47 & 0.55 & 0.63 \\
    \hej{8} & 0.48 & 0.47 & 0.54 & 0.62 \\
    \hej{10} & 0.48 & 0.47 & 0.55 & 0.63 \\
    \midrule
    transforms and plans & 2.24 & 2.74 & 3.08 & 3.36 \\
    one solve & 1.04 & 0.54 & 1.42 & 2.94 \\
    \bottomrule
\end{tabular}

\end{table}


\section{Conclusion}
\label{section:conclusion}
For the lattice stencils as for the regularised kernels, the kernel used so far by Hockney--Eastwood solvers with two unbounded and one periodic direction, the Spietz kernel, is an approximation. Sampling the symbol of the difference operator at the frequencies of the doubled grid amounts to convolving each non-zero periodic mode with its periodised screened Green's function. The source thus acquires periodic images, whose field is an error that no mesh refinement reduces. While the periodic direction is not much longer than the unbounded ones, the image error is of order $\exp(-2\pi d/L_z)$, with $d$ the distance to the nearest image. Beyond, it grows like $(L_z/L_x)^2$ times the net source of each mode. When that net source vanishes, as for a pair of opposite-signed compact sources, it levels off at a value set by the dipole moment instead. The image error is small when the periodic direction is shorter than, or comparable to, the unbounded ones and the source is compact and well centred, though already above the discretisation error of the high-order kernels. It exceeds the solution itself once the periodic direction is several times longer than the unbounded ones, the configuration of vortex tubes, wakes and jets. The usual verification is blind to it: a compact manufactured solution is also the periodic solution, and the discrete residual cannot tell the two kernels apart.

The approximation can be removed at a modest, one-off cost. Every non-zero periodic mode is a screened two-dimensional problem. Its lattice Green's function is obtained in closed form along one unbounded direction, from the one-unbounded LGFs of~\cite{gabbard_lattice_2023}, and by Gauss--Legendre quadrature on a half period along the second. The classical Hockney--Eastwood recipe applied to these slices is exact. The slices being evaluated in physical space, they are not tied to that algorithm either: they can serve any convolution method, such as the lattice fast multipole method by which Hou and Colonius apply the second-order kernel~\cite{hou_adaptive_2024}. The same generator, implemented in \flups{}, delivers the kernels of the dimension-split stencils of order two to eight and of the fourth- and sixth-order Mehrstellen operators. With the regularised one-unbounded kernel as integrand, it also delivers the regularised kernels of order two to ten. Of these only the second-order LGF was previously available~\cite{hou_fast_2025} in its exact form. The kernel is set up once per grid, at the cost of about one solve at most on the grids considered here, and the online solve itself is unchanged.

Exactness has been verified in three ways. The lattice kernels match independent quadratures and satisfy the discrete identity $\Lop G = \Rop\delta$ to $10^{-14}$. A padding test with a white-noise source is at machine precision on every periodic mode, for every kernel, on domains up to eight times longer than wide. The solution of an axially modulated tube source converges at the order of each stencil, up to $2.7\times10^{11}$ unknowns, whereas the Spietz kernel stagnates at the image error whatever the order. On the same source, that error follows the predicted field of the images to $2\%$ or better up to $L_z/L_x = 64$: it grows far beyond the solution itself for a single tube and levels off at the dipole plateau for a pair.

Two directions remain. The lattice kernels are restricted to isotropic lattices by their closed forms along the first unbounded direction, and the mode $k_z = 0$ still relies on the tabulated two-dimensional LGF and its far-field expansion. The quadrature route extends to anisotropic spacings and to the mode $k_z = 0$ itself, with the value at the origin subtracted, and would remove the last precomputed table of the solver. Both are available at second order~\cite{hou_fast_2025, hou_adaptive_2024} and await their transposition to the high-order and compact stencils. All developments described in this paper are publicly available in the open-source \flups{} library~\cite{flups}.

\appendix

\section{Numerical details}
\label{app:numerics}
The closed forms of Step~1 are evaluated as prescribed in~\cite{gabbard_lattice_2023}: the symbols are written as polynomials in $s^2(k) = 4\sin^2(k/2)$ rather than in $\cos k$, and the dimension-split kernels are taken with the small-wavenumber branches of their Appendix~B, which restore full precision when the root nearest the unit circle approaches $1$. In the same spirit, $\sqrt{e_0^2-e_1^2}$ is formed as $\sqrt{\sigL(0;k_y,k_z)\,\sigL(\pi;k_y,k_z)}$: the difference of squares costs a relative error $\epsilon\,e_0^2/(e_0^2-e_1^2)$, which reaches $10^{-10}$ as $(k_y,k_z)\to(0,0)$ on the grids considered below.

One point is not covered there. For the Mehrstellen operators the root reads $\rho = \lambda - \sqrt{\lambda^2-1}$ with $\lambda = e_0/e_1$, and $\lambda\to1^+$ as $(k_y,k_z)\to(0,0)$: forming $\lambda$ from the $\mathcal{O}(1)$ coefficients and then $\lambda - 1$ costs a relative error $\epsilon/(\lambda-1)$ on the difference, of which $\rho - 1\simeq-\sqrt{2(\lambda-1)}$ inherits the square root, $\epsilon/\sqrt{2(\lambda-1)}$, and $\rho^{\,n}$ multiplies it by $n$. Gabbard and van Rees~\cite{gabbard_lattice_2023} note this cancellation and leave it, since in their setting both spectral wavenumbers are frequencies of one grid, so that $\lambda - 1\gtrsim\pi^2/N^2$ and the loss is of order $N\epsilon$ and confined to one mode; their small-wavenumber branch is written for the dimension-split stencils only. Neither argument holds here. The wavenumber $k_y$ is a Gauss--Legendre node, the first of which lies at $\mathcal{O}(1/Q^2)$, so that on the smallest mode of the periodic direction $\lambda - 1\simeq s_z^2/2\simeq 2\pi^2/N_z^2$ is set by $k_z$ alone, and for $N_z = 4096$ and $n\sim N$ the loss reaches $10^{-10}$; and since the nodes $k_y\to0$ carry the long-range part of the outer integral \eqqref{eq:outer}, the error is not local but spreads over the whole slice, above the truncation error of the sixth-order stencils at that resolution.

We therefore never form $\lambda$: the difference $\delta = \lambda - 1 = \sigL(0;k_y,k_z)/e_1$ is evaluated from the slice symbol at $k_x = 0$, a polynomial in $s_y^2$ and $s_z^2$ without constant term ($\sigL(0;k_y,k_z) = s_y^2+s_z^2-s_y^2s_z^2/6$ for both operators), and, with $\tau = \delta+\sqrt{\delta(\delta+2)}$, the root is obtained as $\rho - 1 = -2\delta/\tau$ and $\rho = 1/(1+\tau)$, with $\rho^{\,n} = \exp\bigl(n\log(1+(\rho-1))\bigr)$ through \texttt{log1p} and $\rho^2-1 = (\rho-1)(\rho+1)$: every step is free of cancellation for all $\delta\ge0$. Recovering $\rho$ as $1 + (\rho-1)$ would cancel in the opposite limit $\delta\to\infty$, reached by the fourth-order Mehrstellen operator where $e_1\to0$ and handled in~\cite{gabbard_lattice_2023} by a large-$\lambda$ series; the form $\rho = 1/(1+\tau)$ covers both limits. For that operator $e_1$ changes sign on the curve $\cos k_y + \cos k_z = -1$ and vanishes on it; the region beyond, where $e_1<0$, covers $18\%$ of the $(k_y,k_z)$ plane. There one has $\lambda\le-3$, $\rho = -1/\bigl(|\lambda|+\sqrt{\lambda^2-1}\bigr)$ is well conditioned as it stands, and $\rho<0$ produces the alternation one would obtain by the shift $k\mapsto k+\pi$; on the curve, $e_1 = 0$ gives $\rho = 0$ and $I(n) = \delta_{n0}/e_0$. Written in terms of $(e_0,e_1)$ rather than of $\lambda$, \eqqref{eq:residue} is branch-free across it.

\paragraph{Regularised kernels}
The one-unbounded regularised kernel of~\cite{juul_spietz_regularization_2018} contains the products $\eul^{\pm kx}\,\mathrm{erfc}(\cdot)$; forming the complementary error function as $1-\mathrm{erf}$ loses all digits once its argument exceeds about $4$, which happens on the slices of large $k_z$.

\paragraph{Gauss--Legendre nodes}
The rule \eqqref{eq:Q} is a single panel of up to several thousand nodes. The textbook generator, Newton's method on $P_Q(x)$ with the three-term recurrence and the weights $2/\bigl((1-x^2)P_Q'(x)^2\bigr)$, loses accuracy on the weights next to the endpoints through the cancellation in $1-x^2$~\cite[\S3.1]{hale_fast_2013}: we measure a relative $\epsilon/(1-|x|)\sim\epsilon Q^2$, $10^{-9}$ at $Q = 8192$, which is where the integrand of \eqqref{eq:outer} peaks. As in~\cite[\S3.1]{hale_fast_2013}, we run Newton's method in $\theta$, $x = \cos\theta$, and take the weights from $\mathrm{d}P_Q(\cos\theta)/\mathrm{d}\theta$. Working in $\theta$ is not sufficient with the recurrence, which still rounds $x = \cos\theta$ next to $1$; we therefore evaluate it in the spirit of Reinsch's modification of the Goertzel algorithm~\cite[\S2.3.3]{stoer_introduction_2002}, in the variable $y = 1-x = 2\sin^2(\theta/2)$ and on the increments $u_j = P_j - P_{j-1}$, $u_{j+1} = \bigl(j\,u_j - (2j+1)\,y\,P_j\bigr)/(j+1)$, so that $x$ itself never appears. The nodes $k = \pi\sin^2(\theta/2)$ and the weights $\pi/\bigl(\sin\theta\,P_Q'(\cos\theta)\bigr)^2$ are then accurate to a few ulps at any order, and the rule is built on the accelerator, one thread per pair of nodes.

\paragraph{Truncation of the slices}
Entries of a slice below $10^{-17}$ of its maximum are skipped, the decay rate of every node along $n_x$ being known ($\rho$ for the lattice kernels, $\eul^{-kx}$ for the regularised ones): the kernel carries an absolute floor of $10^{-17}$, of the class of the round-off of the transforms that follow it, which is the sense in which the slices are exact.

\paragraph{Reference solution of the tube source}
The solution \eqqref{eq:wavy_sol} is closed form in $z$ and in its Bessel structure in $r$; the two radial integrals are evaluated by Gauss--Legendre quadrature on geometric panels in $v = t/(1-t)$, $t = (r/R)^2$, in which the bump becomes a plain exponential and the logarithm of $K_0$ sits at a panel endpoint. Against a refined rule the quadrature deviates by less than $4\times10^{-16}$ in absolute value on every mode used, and the circulation of $\omega$ is recovered to $2\times10^{-15}$.

\paragraph{Convergence of the outer quadrature}
\Tref{tab:quadrature} documents the convergence of the Gauss--Legendre rule of \eqqref{eq:outer} and the origin of the second term in \eqqref{eq:Q}: at fixed size $N = 16$, the number of nodes that reaches machine precision on the lowest mode grows with the period $N_z$ once $N_z$ exceeds $2N$, because the boundary layer of $I(n_x;k,k_z)$ at $k = 0$ narrows with $k_z$. Each fourfold increase of $N_z$ doubles it, the $\sqrt{|k_z|}$ scaling of Step~2. With \eqqref{eq:Q} (last column) every entry of the table is at machine precision. With $Q = N$ nodes instead (below the range of the table), the error is $10^{-3}$, and $10^{-10}$ with $Q = 2N$ on a cubic domain.

\begin{table}[htbp]
    \centering
    \caption{Convergence of the single-panel Gauss--Legendre rule of \eqqref{eq:outer} on the \lgf{2} slice $k_z = 2\pi/N_z$ for $N = 16$: maximum error against a $Q = 8192$ reference, for fixed $Q$ and for the rule \eqqref{eq:Q} (its number of nodes in parentheses).}
    \label{tab:quadrature}
    \footnotesize
\setlength{\tabcolsep}{3pt}
\begin{tabular}{r r r r r r r r}
    \toprule
    $N_z$ & $k_z$ & $Q=64$ & $128$ & $256$ & $512$ & $1024$ & \eqqref{eq:Q} \\
    \midrule
       16 & $0.39$ & $2.4\cdot10^{-16}$ & $1.9\cdot10^{-16}$ & $1.6\cdot10^{-16}$ & $1.2\cdot10^{-16}$ & $2.0\cdot10^{-16}$ & $1.5\cdot10^{-16}$ ($96$) \\
       64 & $9.8\cdot10^{-2}$ & $3.2\cdot10^{-15}$ & $2.5\cdot10^{-16}$ & $2.0\cdot10^{-16}$ & $2.0\cdot10^{-16}$ & $3.9\cdot10^{-16}$ & $1.9\cdot10^{-16}$ ($160$) \\
      256 & $2.5\cdot10^{-2}$ & $1.7\cdot10^{-8}$ & $1.7\cdot10^{-15}$ & $3.0\cdot10^{-16}$ & $3.6\cdot10^{-16}$ & $5.3\cdot10^{-16}$ & $4.7\cdot10^{-16}$ ($288$) \\
     1024 & $6.1\cdot10^{-3}$ & $8.7\cdot10^{-5}$ & $1.6\cdot10^{-8}$ & $2.0\cdot10^{-15}$ & $2.3\cdot10^{-16}$ & $5.8\cdot10^{-16}$ & $2.6\cdot10^{-16}$ ($543$) \\
     4096 & $1.5\cdot10^{-3}$ & $7.6\cdot10^{-3}$ & $9.1\cdot10^{-5}$ & $1.6\cdot10^{-8}$ & $2.0\cdot10^{-15}$ & $3.0\cdot10^{-16}$ & $2.6\cdot10^{-16}$ ($1054$) \\
    16384 & $3.8\cdot10^{-4}$ & $5.8\cdot10^{-2}$ & $7.7\cdot10^{-3}$ & $9.3\cdot10^{-5}$ & $1.6\cdot10^{-8}$ & $2.5\cdot10^{-15}$ & $7.4\cdot10^{-16}$ ($2075$) \\
    \bottomrule
\end{tabular}

\end{table}

\section*{Acknowledgments}
    Generative AI tools were used in this work under the direction and review of the authors. The authors assume responsibility for all content.
    Computational resources have been provided by the supercomputing facilities of the Université catholique de Louvain (CISM/UCL) and the Consortium des Équipements de Calcul Intensif en Fédération Wallonie Bruxelles (CÉCI) funded by the Fond de la Recherche Scientifique de Belgique (F.R.S.-FNRS) under convention 2.5020.11 and by the Walloon Region.
    The present research benefited from computational resources made available on Lucia, the Tier-1 supercomputer of the Walloon Region, infrastructure funded by the Walloon Region under the grant agreement n°1910247.
    We acknowledge LUMI-BE\footnote{LUMI-BE is joint effort from BELSPO (federal), SPW Économie, Emploi, Recherche (Wallonia), Department of Economy, Science \& Innovation (Flanders) and Innoviris (Brussels).} for awarding this project access to the LUMI supercomputer, owned by the EuroHPC Joint Undertaking, hosted by CSC (Finland) and the LUMI consortium through a LUMI-BE Regular Access call.

\bibliographystyle{siam}
\bibliography{AllPublications}

\end{document}